\magnification\magstep1
\input picmac
\input epsf
\def\mpfig#1{\hbox{$\vcenter{\epsfbox{sand.#1}}$}}
% This file uses PostScript illustrations sand.1 thru sand.5 made by MetaPost

\font\rhfont=cmcsc8
\headline=
 {\ifnum\pageno>1
   \rhfont The Electronic Journal of Combinatorics 1 (1994), \# A1
   \hfil \tenrm\folio
  \else \hfil \fi}
\footline={\hfil}

% Reference numbers:
\def\Ala{1} % Alizadeh 1991
\def\Alb{2} % Alizadeh 1993
\def\ALMSS{3} % Arora et al
\def\Bol{4} % Bollobas
\def\GLSell{5} % Groetschel et al 1981
\def\GLSrelax{6} % Groetschel et al 1986
\def\GLSbook{7} % Groetschel et al 1988
\def\House{8} % Householder
\def\Ju{9} % Juhasz
\def\KK{10} % Kashin
\def\Ko{11} %Koniagin
\def\Lov{12} % Lovasz 1979
\def\LL{13} % Lovasz 1986
\def\LLL{14} % Lovasz to appear
\def\LS{15} % Lovasz, Schrijver
\def\Ov{16} % Overton
\def\Padberg{17}
\def\Shannon{18}

\let\SQRT=\sqrt
\def\sqrt#1{\SQRT{\vphantom t#1}}
\catcode`\@=\active
\def@{\ifmmode\mskip1mu\else\char`\@ \fi}
\def\0{^{\vphantom T}}

\parskip3pt
\parindent20pt
\baselineskip14pt

\def\rar{\rightarrow}
\def\meno{\medskip\noindent}
\def\adj{\mathrel-\joinrel\joinrel\mathrel-}   %long minus
\def\noadj{\hskip3pt\not\hskip-3pt\adj}  %not long minus
\def\ones{{1\mskip-4.5mu\rm l}}
\def\bib{\par\noindent\hangindent 25pt}
\def\disleft#1:#2:#3\par{\par\hangindent#1\noindent
	 \hbox to #1{#2 \hfill \hskip .1em}\ignorespaces#3\par}
\def\pfbox
  {\hbox{\hskip 3pt\lower1pt\vbox{\hrule
  \hbox to 7pt{\vrule height 7pt\hfill\vrule}
  \hrule}}\hskip3pt}
\def\dddots{\mathinner{\mskip1mu\raise1pt\vbox{\kern7pt\hbox{.}}\mskip2mu
    \raise4pt\hbox{.}\mskip2mu\raise7pt\hbox{.}\mskip1mu}}
\def\ddddots{\mathinner{\mskip1mu\raise7pt\vbox{\kern7pt\hbox{.}}\mskip2mu
    \raise4pt\hbox{.}\mskip2mu\raise1pt\hbox{.}\mskip1mu}}

\pageno=0

\centerline{\bf The Sandwich Theorem}
\bigskip
\centerline{Donald E. Knuth}

\bigskip\bigskip\bigskip\bigskip
{\narrower\narrower\narrower\smallskip\noindent
{\bf Abstract:}\enspace
This report contains expository notes about a function $\vartheta(G)$ that is
popularly known as the Lov\'asz number of a graph~$G$. There are many ways to
define $\vartheta(G)$, and the surprising variety of different
characterizations indicates in itself  that $\vartheta(G)$
should be interesting. But the most interesting property of $\vartheta(G)$ is
probably the fact that it can be computed efficiently, although it lies
``sandwiched'' between other classic graph numbers whose computation is
NP-hard. I~have tried to make these notes self-contained so that they might
serve as an elementary introduction to the growing literature on Lov\'asz's
fascinating function.
\smallskip}

\vfill\eject

\line{{\bf The Sandwich Theorem}
      \hfill DEK notes last revised 6 December 1993}

\medskip
\begingroup \advance\baselineskip 0pt minus1pt
\halign{\hfil#\quad&\hbox to .9\hsize{# \dotfill\quad}&\hfil#\cr
0.&Preliminaries&2\cr
1.&Orthogonal labelings&4\cr
2.&Convex labelings&5\cr
3.&Monotonicity&6\cr
4.&The theta function&6\cr
5.&Alternative definitions of $\vartheta$&8\cr
6.&Characterization via eigenvalues&8\cr
7.&A complementary characterization&9\cr
8.&Elementary facts about cones&11\cr
9.&Definite proof of a semidefinite fact&12\cr
10.&Another characterization&13\cr
11.&The final link&14\cr
12.&The main theorem&14\cr
13.&The main converse&15\cr
14.&Another look at {\tt TH}&16\cr
15.&Zero weights&17\cr
16.&Nonzero weights&17\cr
17.&Simple examples&19\cr
18.&The direct sum of graphs&19\cr
19.&The direct cosum of graphs&20\cr
20.&A direct product of graphs&21\cr
21.&A direct coproduct of graphs&22\cr
22.&Odd cycles&24\cr
23.&Comments on the previous example&27\cr
24.&Regular graphs&27\cr
25.&Automorphisms&28\cr
26.&Consequence for eigenvalues&30\cr
27.&Further examples of symmetric graphs&30\cr
28.&A bound on $\vartheta$&31\cr
29.&Compatible matrices&32\cr
30.&Antiblockers&35\cr
31.&Perfect graphs&36\cr
32.&A characterization of perfection&37\cr
33.&Another definition of $\vartheta$&39\cr
34.&Facets of {\tt TH}&40\cr
35.&Orthogonal labelings in a perfect graph&42\cr
36.&The smallest non-perfect graph&43\cr
37.&Perplexing questions&45\cr}

\endgroup
\vfill\eject

\centerline{\bf The Sandwich Theorem}

\bigskip\noindent
It is NP-complete to compute $\omega(G)$, the size of the largest
clique in a graph~$G$, and it is NP-complete to compute $\chi(G)$, the
minimum number of colors needed to color the vertices of~$G$. But 
Gr\"otschel,  Lov\'asz, and Schrijver proved
[\GLSell]
that we can compute in polynomial time a real number that is ``sandwiched''
between these hard-to-compute integers:
$$\omega(G)\leq\vartheta(\overline{G})\leq\chi(G)\,.\eqno(\ast)$$
Lov\'asz [\LL] called this a ``sandwich theorem.'' The book [\GLSbook]
develops further facts about the function $\vartheta(G)$ and shows that
it possesses many interesting properties. Therefore I~think it's
worthwhile to study $\vartheta(G)$ closely, in hopes of getting 
acquainted with it and finding faster
ways to compute~it.

Caution: The function called $\vartheta(G)$ in [\LL] is called
$\vartheta(\overline{G})$ in [\GLSbook] and~[\Lov]. I~am following the latter
convention because it is more likely to be adopted by other
researchers---[\GLSbook] is a classic book that contains complete proofs,
while [\LL] is simply an extended abstract.

In these notes I am mostly following [\GLSbook] and [\Lov] with minor
simplifications and a few additions. I~mention several natural
problems that I~was not able to solve immediately although I~expect
(and fondly hope) that they will be resolved before I~get to writing
this portion of my forthcoming book on Combinatorial Algorithms. 
I'm grateful to many people---especially to Martin Gr\"otschel and L\'aszl\'o
Lov\'asz---for their comments on my first drafts of this material.

These notes are in
numbered sections, and there is at most one Lemma, Theorem, Corollary,
or Example in each section. Thus, ``Lemma~2'' will mean ``the lemma in
section~2''.

\meno
{\bf 0. Preliminaries.}\quad
Let's begin slowly by defining some notational conventions and by
stating some
basic things that will be assumed without proof. All vectors in these notes
will be regarded as column vectors, indexed either by the vertices of a graph
or by integers. 
The notation $x\geq y$, when $x$ and~$y$ are vectors, will mean
that $x_v\geq y_v$ for all~$v$. If $A$ is a matrix, $A_v$~will denote
column~$v$, and $A_{uv}$ will be the element in row~$u$ of column~$v$. The
zero vector and the zero matrix and zero itself will all be denoted by~0. 

We will use several properties of 
matrices and vectors of real numbers that are
familiar to everyone who works with linear algebra but not to everyone who
studies graph theory, so it seems wise to list them here:

\medskip
(i)\quad The {\it dot product\/} of (column) vectors~$a$ and~$b$ is
$$a\cdot b=a^Tb\,;\eqno(0.1)$$
the vectors are {\it orthogonal\/} (also called perpendicular) 
if $a\cdot b=0$.
The {\it length\/} of vector~$a$ is
$$\Vert a\Vert=\sqrt{a\cdot a}\,.\eqno(0.2)$$
Cauchy's inequality asserts that
$$a\cdot b\leq \Vert a\Vert\;\Vert b\Vert\,;\eqno(0.3)$$
equality holds iff $a$ is a scalar multiple of $b$ 
or $b=0$. Notice that if $A$
is any matrix we have
$$(A^T\!A)_{uv}=\sum_{k=1}^n(A^T)_{uk}A_{kv}=\sum_{k=1}^nA_{ku}A_{kv}=A_u\cdot
A_v\,;\eqno(0.4)$$
in other words, the elements of $A^T\!A$ represent all dot products of the
columns of~$A$.

\medskip
(ii)\quad
An {\it orthogonal matrix\/} is a square
matrix~$Q$ such that $Q^TQ$ is the identity
matrix~$I$. Thus, by (0.4), $Q$~is orthogonal
iff its columns are unit vectors
perpendicular to each other. The transpose of an orthogonal matrix is
orthogonal, because the condition $Q^TQ=I$ implies that $Q^T$ is the inverse
of~$Q$, hence $QQ^T=I$. 

\medskip
(iii)\quad
A given matrix $A$ is {\it symmetric\/} (i.e., $A=A^T$) iff it can be
expressed in the form
$$A=QDQ^T\eqno(0.5)$$
where $Q$ is orthogonal and $D$ is a diagonal matrix. Notice that (0.5) is
equivalent to the matrix equation
$$AQ=QD\,,\eqno(0.6)$$
which is equivalent to the equations
$$AQ_v=Q_v\lambda_v$$
for all $v$, where $\lambda_v=D_{vv}$. Hence the diagonal elements of~$D$ are
the eigenvalues of~$A$ and the columns of~$Q$ are the corresponding
eigenvectors. 

Properties (i), (ii), and (iii) 
are proved in any textbook of linear algebra. We
can get some practice using these concepts by giving a constructive proof of
another well known fact:

\proclaim
Lemma. Given $k$ mutually perpendicular unit vectors, there is an orthogonal
matrix having these vectors as the first $k$~columns.

\noindent{\bf Proof.}\quad 
Suppose first that $k=1$ and that $x$ is a $d$-dimensional vector with $\Vert
x\Vert=1$. If $x_1=1$ we have $x_2=\cdots =x_d=0$, so the orthogonal matrix
$Q=I$ satisfies the desired condition. Otherwise we let
$$y_1=\sqrt{(1-x_1)/2}\,,\qquad y_j=-x_j/(2y_1)\quad{\rm for}\ 1<j\leq d\,.
\eqno(0.7)$$ 
Then
$$y^T\!y=\Vert y\Vert^2=y_1^2+{x_2^2+\cdots +x_d^2\over 4y_1^2}={1-x_1\over
2}+{1-x_1^2\over 2(1-x_1)}=1\,.$$
And $x$ is the first column of the Householder [\House] matrix
$$Q=I-2yy^T\,,\eqno(0.8)$$
which is easily seen to be orthogonal because
$$Q^TQ=Q^2=I-4yy^T+4yy^Tyy^T=I\,.$$

Now suppose the lemma has been proved for some $k\geq 1$; we will show
how to increase $k$ by~1. 
Let $Q$ be an orthogonal matrix and let $x$ be a unit
vector perpendicular to its first $k$~columns. We want to construct an
orthogonal matrix~$Q'$ agreeing with $Q$ in columns~1 to~$k$ and having $x$ in
column $k+1$. Notice that
$$Q^Tx=\pmatrix{0\cr \vdots\cr 0\cr y\cr}$$
by hypothesis, where there are 0s in the first $k$ rows. The
$(d-k)$-dimensional vector~$y$ has squared length
$$\Vert y\Vert^2=Q^Tx\cdot Q^Tx=x^TQQ^Tx=x^Tx=1\,,$$
so it is a unit vector. (In particular, $y\neq 0$, so we must have $k<d$.)
Using the construction above, 
we can find a $(d-k)\times (d-k)$ orthogonal matrix~$R$
with $y$ in its first column. Then the matrix
$$Q'=Q\pmatrix{1\cr &\ddddots&&0\cr &&1\cr &0&&R\cr}$$
does what we want.\ \pfbox

\meno
{\bf 1. Orthogonal labelings.}\quad
Let $G$ be a graph on the vertices $V$. If $u$ and~$v$ are distinct elements
of~$V$, the notation $u\adj v$ means that
they are adjacent in~$G$; $u\noadj v$ means they are not.

An assignment of vectors $a_v$ to each vertex $v$ is called an {\it orthogonal
labeling\/} of~$G$ if $a_u\cdot a_v=0$ whenever $u\noadj v$. In other words,
whenever $a_u$ is not perpendicular to~$a_v$ in the labeling, 
we must have $u\adj v$ in the graph. The
vectors  may have any desired
dimension~$d$; the components of~$a_v$ are~$a_{jv}$ for $1\leq j\leq
d$. Example: $a_v=0$ for all~$v$ always works trivially.

The {\it cost\/} $c(a_v)$ of a vector $a_v$ 
 in an orthogonal labeling is defined to
be~0 if $a_v=0$, otherwise
$$c(a_v)={a_{1v}^2\over\Vert a_v\Vert^2}={a_{1v}^2\over
a^2_{1v}+\cdots +a^2_{dv}}\,.$$
Notice that we can multiply any vector $a_v$ by a nonzero scalar~$t_v$
without changing its cost, and without violating the orthogonal
labeling property. 
We can also get rid of a zero vector by increasing $d$ by~1 and adding a new
component~0 to each vector, except that the zero vector gets the new
component~1. 
In particular, we can if we like assume that all
vectors have unit length. Then the cost will be $a^2_{1v}$.

\proclaim
Lemma.  If $S\subseteq V$ is a {\it stable set\/} of vertices (i.e., no
two vertices of~$S$ are adjacent) and if $a$ is an orthogonal labeling
then
$$\sum_{v\in S}c(a_v)\leq 1\,.\eqno(1.1)$$

\noindent
{\bf Proof.}\quad
We can assume that $\Vert a_v\Vert =1$ for all $v$. Then the vectors~$a_v$ for
$v\in S$ must be mutually orthogonal, and Lemma~0 tells us we can find a
$d\times d$ orthogonal matrix~$Q$ with these vectors as its leftmost columns. 
The sum of the
costs will then be at most $q^2_{11}+q^2_{12}+\cdots +q^2_{1d}=1$. \pfbox 

\medskip\noindent
Relation (1.1) makes it
possible for us to study stable sets geometrically.

\meno
{\bf 2. Convex labelings.}\quad
An assignment $x$ of real numbers $x_v$ to the vertices~$v$ of~$G$ is
called a {\it real labeling\/} of~$G$.
 Several families of such labelings will be of importance to us:

\disleft40pt::
The {\it characteristic labeling\/} for $U\subseteq V$ has
$x_v=\cases{1&if $v\in U$;\cr
0&if $v\notin U$.\cr}$

\disleft40pt::
A {\it stable labeling\/} is a characteristic labeling for a stable
set.

\disleft40pt::
A {\it clique labeling\/} is a characteristic labeling for a clique (a~set of
mutually adjacent 

\disleft60pt::
vertices).

\disleft40pt::
{\tt STAB}$(G)$ is the smallest convex set containing all stable labelings,

\disleft60pt::
i.e., {\tt STAB}$(G)=$ convex hull $\{\,x\mid x$ is a stable
labeling of~$G\,\}$.

\disleft40pt::
${\tt QSTAB}(G) =\{\,x\geq 0
\mid\sum_{v\in Q}x_v\leq 1\hbox{ for all cliques $Q$ of $G$}\,\}\,.$

\disleft40pt::
${\tt TH}(G)=\{\,x\geq 0\mid \sum_{v\in V}c(a_v)x_v \leq 1\hbox{ for
all orthogonal labelings $a$ of~$G$}\,\}\,.$

\proclaim
Lemma. {\tt TH} is sandwiched between {\tt STAB} and {\tt QSTAB}:
$${\tt STAB}(G)\subseteq {\tt TH}(G)\subseteq{\tt
QSTAB}(G)\,.\eqno(2.1)$$

\noindent
{\bf Proof.}\quad
Relation (1.1) tells that every stable labeling belongs to {\tt
TH}$(G)$. Since {\tt TH}$(G)$ is obviously convex, it must contain the
convex hull {\tt STAB}$(G)$. On the other hand, every clique labeling is an
orthogonal labeling of dimension~1.
 Therefore every constraint of {\tt QSTAB}$(G)$
is one of the constraints of {\tt TH}$(G)$. \ \pfbox

\meno
{\bf Note:}\quad
{\tt QSTAB} first defined by Shannon [\Shannon], 
 and the first systematic study of {\tt STAB} was undertaken by Padberg
[\Padberg]. {\tt TH} was first defined by Gr\"otschel, Lov\'asz, and Schrijver
in~[\GLSrelax].

\meno
{\bf 3. Monotonicity.}\quad
Suppose $G$ and $G'$ are graphs on the same vertex set~$V$, with
$G\subseteq G'$ (i.e., $u\adj v$ in~$G$ implies $u\adj v$ in~$G'$).
Then 

\disleft40pt::
every stable set in~$G'$ is stable in~$G$, hence ${\tt
STAB}(G)\supseteq {\tt STAB}(G')$;

\disleft40pt::
every clique in $G$ is a clique in $G'$, hence ${\tt
QSTAB}(G)\supseteq {\tt QSTAB}(G')$;

\disleft40pt::
every orthogonal labeling of $G$ is an orthogonal labeling of $G'$,

\disleft60pt::
hence ${\tt TH}(G)\supseteq {\tt TH}(G')$.

\noindent
In particular, if $G$ is the empty graph $\overline{K}_n$ on $\vert V\vert
=n$ vertices, all sets are stable and all cliques have size $\leq 1$,
hence

\disleft40pt::
${\tt STAB}(\overline{K}_n)={\tt TH}(\overline{K}_n)
={\tt QSTAB}(\overline{K}_n)
=\{\,x\mid 0\leq x_v\leq 1\hbox{ for all }v\,\}$, the $n$-cube.

\noindent
If $G$ is the complete graph $K_n$, all stable sets have size $\leq
1$ and there is an $n$-clique, so

\disleft40pt::
${\tt STAB}(K_n)={\tt TH}(K_n)={\tt QSTAB}(K_n)=\left\{\,x\geq 0\mid 
\sum_{v}x_v\leq 1\,\right\}$, the $n$-simplex.

\noindent
Thus all the convex sets {\tt STAB}$(G)$, {\tt TH}$(G)$, {\tt
QSTAB}$(G)$ lie between the $n$-simplex and the $n$-cube.

Consider, for example, the case $n=3$. Then there are three coordinates, so we
can visualize the sets in 3-space (although there aren't many interesting
graphs). The {\tt QSTAB} of
{\unitlength=10pt\beginpicture(3,.5)(-.5,-.25)
\def\putdisk(#1,#2){\put(#1,#2){\disk{.4}}}%
\def\putlab(#1,#2)#3{\put(#1,#2){\makebox(0,0){$\scriptstyle\smash{#3}$}}}%
\putdisk(0,0)\putdisk(1,0)\putdisk(2,0)
\sevenrm\putlab(0,.5)x\putlab(1,.5)y\putlab(2,.5)z%
\put(0,0){\line(1,0){2}}\endpicture}
is obtained from the unit cube by
restricting the coordinates to $x+y\leq 1$ and $y+z\leq 1$; we can think of
making two cuts in a piece of cheese:
$$
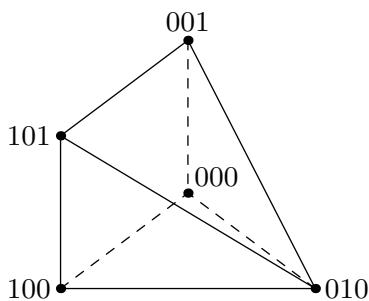
\mpfig5
$$
The vertices $\{000, 100, 010, 001, 101\}$ 
correspond to the stable labelings,
so once again we have ${\tt STAB}(G)={\tt TH}(G)={\tt QSTAB}(G)$.

\meno
{\bf 4. The theta function.}\quad
The function $\vartheta(G)$ mentioned in the introduction is a special
case of a two-parameter function $\vartheta(G,w)$, where $w$ is a
nonnegative real labeling:
$$\eqalignno{\vartheta(G,w)&=\max\{\,w\cdot x\mid x\in {\tt
TH}(G)\,\}\,;&(4.1)\cr
\noalign{\smallskip}
\vartheta(G)&=\vartheta(G,\ones)\,\ \hbox{where $\ones$ is the
labeling $w_v=1$ for all $v$.}&(4.2)\cr}$$
This function, called the {\it Lov\'asz number\/} of~$G$ (or the {\it weighted
Lov\'asz number\/} when $w\neq \ones$), tells us about 1-dimensional
projections of the $n$-dimensional convex set ${\tt TH}(G)$.

Notice, for example, that the monotonicity  properties of \S3 tell us
$$G\subseteq G'\ \Rightarrow\ \vartheta(G,w)\geq
\vartheta(G',w)\eqno(4.3)$$ for all $w\geq 0$. It is also obvious that
$\vartheta$ is monotone in its other parameter:
$$w\leq w'\ \Rightarrow\ \vartheta(G,w)\leq
\vartheta(G,w')\,.\eqno(4.4)$$
The smallest possible value of $\vartheta$ is
$$\vartheta(K_n,w)=\max\{w_1,\ldots,w_n\}\,;\quad
\vartheta(K_n)=1\,.\eqno(4.5)$$ 
The largest possible value is
$$\vartheta(\overline{K}_n,w)=w_1+\cdots +w_n\,;\quad
\vartheta(\overline{K}_n)=n\,.\eqno(4.6)$$

Similar definitions can be given for {\tt STAB} and {\tt QSTAB}:
$$\eqalignno{\alpha(G,w)&=\max\{\,w\cdot x\mid x\in{\tt STAB}(G)\,\}\,,
\qquad\hskip4pt \alpha(G)=\alpha(G,\ones)\,;&(4.7)\cr
\noalign{\smallskip}
\kappa(G,w)&=\max\{\,w\cdot x\mid x\in {\tt QSTAB}(G)\,\}\,,
\qquad \kappa(G)=\kappa(G,\ones)\,.&(4.8)\cr}$$
Clearly $\alpha(G)$ is the size of the largest stable set in~$G$, because 
every
stable labeling~$x$ corresponds to a stable set with $\ones\cdot x$ vertices.
It is also easy to see that $\kappa(G)$ is at most $\overline{\chi}(G)$, the
smallest number of cliques that cover the vertices of~$G$. For if the vertices 
can be partitioned into $k$~cliques $Q_1,\ldots,Q_k$ and if $x\in{\tt
QSTAB}(G)$, we have
$$\ones\cdot x=\sum_{v\in Q_1}x_v+\cdots +\sum_{v\in Q_k}x_v\leq k\,.$$
Sometimes $\kappa(G)$ is less than 
$\overline{\chi}(G)$. For example, consider  the
cyclic graph~$C_n$, with vertices $\{0,1,\ldots,n-1\}$ and $u\adj v$ iff 
$u\equiv v\pm 1$~(mod~1). Adding up the inequalities $x_0+x_1\leq 1,\ldots,
x_{n-2}+x_{n-1}\leq 1,x_{n-1}+x_0\leq 1$ of {\tt QSTAB} gives
$2(x_0+\cdots +x_{n-1})\leq n$, and this upper bound is achieved when
 all~$x$'s are~${1\over 2}$; hence $\kappa(C_n)={n\over 2}$, if $n>3$. But
$\overline{\chi}(G)$ is always an integer, and $\overline{\chi}(C_n)
=\left\lceil{n\over 2}\right\rceil$ is greater than $\kappa(C_n)$ when $n$ is
odd.

Incidentally, these remarks establish the ``sandwich inequality'' $(\ast)$
stated in the introduction, because
$$\alpha(G)\leq\vartheta(G)\leq\kappa(G)\leq\overline{\chi}(G)\eqno(4.9)$$
and $\omega(G)=\alpha(\overline{G})$, $\overline{\chi}(G)=
\chi(\overline{G})$.

\meno
{\bf 5. Alternative definitions of $\vartheta$.}\quad
Four additional functions $\vartheta_1$, $\vartheta_2$, $\vartheta_3$,
$\vartheta_4$ are defined in~[\GLSbook], and they all turn out to be
identical to~$\vartheta$. Thus, we can understand $\vartheta$ in many
different ways; this may help us compute~it.

We will show, following [\GLSbook], that if $w$ is any fixed nonnegative
real labeling of~$G$, the inequalities
$$\vartheta(G,w)\leq \vartheta_1(G,w)\leq\vartheta_2(G,w)\leq
\vartheta_3(G,w)\leq\vartheta_4(G,w)\leq\vartheta(G,w)\eqno(5.1)$$
can be proved. Thus we will establish the theorem of~[\GLSbook], and all
inequalities in our proofs will turn out to be equalities. We will
introduce the alternative definitions~$\vartheta_k$ one at a time; any
one of these definitions could have been taken as the starting point.
First,
$$\vartheta_1(G,w)=\min_a\;\max_v\bigl(w_v/c(a_v)\bigr)\,,\hbox{
over all orthogonal labelings $a$}.\eqno(5.2)$$
Here we regard $w_v/c(a_v)=0$ when $w_v=c(a_v)=0@$; but the max
 is~$\infty$ if some $w_v>0$ has $c(a_v)=0$.

\proclaim
Lemma. $\vartheta(G,w)\leq \vartheta_1(G,w)$.

\noindent
{\bf Proof.}\quad
Suppose $x\in{\tt TH}(G)$ maximizes $w\cdot x$, and suppose $a$ is an
orthogonal labeling that achieves the minimum value
$\vartheta_1(G,w)$. Then
$$\vartheta(G,w)=w\cdot x=\sum_vw_vx_v\leq \biggl(\max_v\,
{w_v\over c(a_v)}\biggr)\,\sum_v c(a_v)x_v\leq
\max_v\,{w_v\over
c(a_v)}=\vartheta_1(G,w)\,.\ \pfbox$$

Incidentally, the fact that all inequalities are exact
 will imply later  that every nonzero weight vector~$w$ has an
orthogonal labeling~$a$ such that
$$c(a_v)={w_v\over\vartheta(G,w)}\quad\hbox{for all $v$}.\eqno(5.3)$$
We will restate such consequences of (5.1) later,
 but it may be helpful to keep that future goal in mind.

\meno
{\bf 6. Characterization via eigenvalues.}\quad
The second variant of~$\vartheta$ is rather different; this is the only
one Lov\'asz chose to mention in~[\LL].

We say that $A$ is a {\it feasible matrix\/} for $G$ and $w$ if $A$ is
indexed by vertices and
$$\eqalignno{%
&\hbox{$A$ is real and symmetric;}\cr
\noalign{\smallskip}
&\hbox{$A_{vv}=w_v$ for all $v\in V$;}\cr
\noalign{\smallskip}
&\hbox{$A_{uv}=\sqrt{w_uw_v}$ whenever $u\noadj v$ in $G$}&(6.1)\cr}$$
The other elements of $A$ are unconstrained (i.e., they can be
anything between $-\infty$ and $+\infty$).

If $A$ is any real, symmetric matrix, let $\Lambda(A)$ be its maximum
eigenvalue. This is well defined because all eigenvalues of~$A$ are
real. Suppose $A$ has eigenvalues $\{\lambda_1,\ldots,\lambda_n\}$; then
$A=Q\,{\rm diag}(\lambda_1,\ldots,\lambda_n)@Q^T$ for some orthogonal~$Q$, and
$\Vert Qx\Vert=\Vert x\Vert$ for all vectors~$x$, so there is a nice way to
characterize $\Lambda(A)$:
$$\Lambda(A)=\max\{\,x^T\!Ax\,\mid\,\Vert x\Vert =1\,\}\,.\eqno(6.2)$$
Notice that $\Lambda(A)$ might not be the largest eigenvalue in
absolute value. We now let
$$\vartheta_2(G,w)=\min\{\,\Lambda(A)\mid A\hbox{ is a feasible matrix
for $G$ and $w$}\,\}\,.\eqno(6.3)$$

\proclaim
Lemma. $\vartheta_1(G,w)\leq \vartheta_2(G,w)$.

\noindent
{\bf Proof.}\quad
Note first that the trace ${\rm tr}\,A=\sum_v w_v\geq 0$ for any
feasible matrix~$A$. The trace is also well-known to be the sum of the
eigenvalues; this fact is an easy consequence of the identity
$${\rm tr}\,XY=\sum_{j=1}^m\;\sum_{k=1}^nx_{jk}y_{kj}={\rm tr}\,YX
\eqno(6.4)$$
valid for any matrices $X$ and $Y$ of respective sizes $m\times n$ and
$n\times m$. In particular, $\vartheta_2(G,w)$ is always $\geq 0$, and
it is $=0$ if and only if $w=0\ 
\bigl($when also $\vartheta_1(G,w)=0\bigr)$.

So suppose $w\neq 0$ and let $A$ be a feasible matrix that attains the
minimum value $\Lambda(A)=\vartheta_2(G,w)=\lambda >0$. Let
$$B=\lambda I-A\,.\eqno(6.5)$$
The eigenvalues of $B$ are $\lambda$ minus the eigenvalues of $A$.
$\bigl($For if 
$A=Q\,{\rm diag}(\lambda_1,\ldots,\lambda_n)@Q^T$  then
$B=Q\,{\rm diag}(\lambda-\lambda_1,\ldots,\lambda-\lambda_n)@Q^T.\bigr)$
Thus they are all nonnegative; such a matrix~$B$ is called {\it positive
semidefinite}. By (0.5) we can write
$$B=X^TX\,,\qquad\hbox{i.e.,}\quad B_{uv}=x_u\cdot x_v\,,\eqno(6.6)$$
when $X={\rm
diag}(\sqrt{\lambda-\lambda_1},\ldots,\sqrt{\lambda-\lambda_n}\,)@Q^T$. 

Let $a_v=(\sqrt{w_v}\,,x_{1v},\ldots,x_{rv})^T$. Then
$c(a_v)=w_v/\,\Vert a_v\Vert^2=w_v/(w_v+x^2_{1v}+\cdots +x^2_{rv})$
and $x^2_{1v}+\cdots +x^2_{rv}=B_{vv}=\lambda-w_v$, hence
$c(a_v)=w_v/\lambda$. Also if $u\noadj v$ we have $a_u\cdot a_v=
\sqrt{w_uw_v}+x_u\cdot x_v =
\sqrt{w_uw_v}+B_{uv}=\sqrt{w_uw_v}-A_{uv}=0$. Therefore $a$ is an
orthogonal labeling and $\max_v\,w_v/c(a_v)=\lambda\geq
\vartheta_1(G,w)$.\ \pfbox

\meno
{\bf 7. A complementary characterization.}\quad
Still another variation is based on orthogonal labelings of the
complementary graph~$\overline{G}$.

In this case we let $b$ be an orthogonal labeling of~$\overline{G}$,
normalized so that $\sum_v\Vert b_v\Vert^2=1$, and we let
$$\eqalignno{\vartheta_3(G,w)=\max\left\{\,\sum_{u,v}\right.
(\sqrt{w_u}\,b_u)&\cdot (\sqrt{w_v}\,b_v)\biggm| \cr
\noalign{\smallskip}
&\left.\phantom{\sum_{u,v}}\hskip-2em
b\hbox{ is a normalized orthogonal labeling of 
$\overline{G}$}\,\right\}.&(7.1)\cr}$$
A normalized orthogonal labeling $b$ is equivalent to an $n\times n$
symmetric positive semidefinite matrix~$B$, where $B_{uv}=b_u\cdot b_v$
is zero when $u\adj v$, and where ${\rm tr}\,B=1$.

\proclaim
Lemma. $\vartheta_2(G,w)\leq\vartheta_3(G,w)$.

This lemma is the ``heart'' of the proof that all $\vartheta$s are
equivalent, according to~[\GLSbook]. It relies on a fact about positive
semidefinite matrices that we will prove in \S9.

\meno
\proclaim
Fact.
If $A$ is a symmetric matrix such that $A\cdot B\geq 0$ for all
symmetric positive semi\-definite~$B$ with $B_{uv}=0$ for $u\adj v$,
then $A=X+Y$ where $X$ is symmetric positive semidefinite and $Y$ is
symmetric and $Y_{vv}=0$ for all~$v$ and $Y_{uv}=0$ for $u\noadj v$.

Here $C\cdot B$ stands for the dot product of matrices, i.e., the sum
$\sum_{u,v}C_{uv}B_{uv}$, which can also be written ${\rm tr}\,C^TB$.
The stated fact is a duality principle for quadratic programming.

Assuming the Fact, let $W$ be the matrix with
$W_{uv}=\sqrt{w_uw_v}$, and let $\vartheta_3=\vartheta_3(G,w)$. By
definition (7.1), if $b$ is any nonzero orthogonal labeling
of~$\overline{G}$ (not necessarily normalized), we have
$$\sum_{u,v}(\sqrt{w_u}\,b_u)\cdot(\sqrt{w_v}\,b_v)\leq\vartheta_3\sum_v
\Vert b_v\Vert^2\,.\eqno(7.2)$$
In matrix terms this says $W\cdot B\leq(\vartheta_3 I)\cdot B$ for all
symmetric positive semidefinite~$B$ with $B_{uv}=0$ for $u\adj v$.
The Fact now tells us we can write
$$\vartheta_3 I-W=X+Y\eqno(7.3)$$
where $X$ is symmetric positive semidefinite, $Y$~is symmetric and
diagonally zero, and $Y_{uv}=0$ when $u\noadj v$. Therefore the
matrix~$A$ defined by
$$A=W+Y=\vartheta_3 I-X$$
is a feasible matrix for $G$, and $\Lambda(A)\leq\vartheta_3$. This
completes the proof that $\vartheta_2(G,w)\leq \vartheta_3(G,w)$,
because $\Lambda(A)$ is an upper bound on $\vartheta_2$ by definition
of~$\vartheta_2$.\ \pfbox

\meno
{\bf 8. Elementary facts about cones.}\quad
A {\it cone\/} in $N$-dimensional space is a set of vectors closed
under addition and under multiplication by nonnegative scalars. (In
particular, it is convex: If $c$ and~$c'$ are in cone~$C$ and
$0<t <1$ then $t c$ and $(1-t)c'$ are in~$C$, hence
$t c+(1-t)c'\in C$.) A~{\it closed cone\/} is a cone that
is also closed under taking limits.

\meno
{\bf F1.}\quad
If $C$ is a closed convex set and $x\notin C$, there is a hyperplane
separating $x$ from~$C$. This means there is a vector~$y$ and a
number~$b$ such that $c\cdot y\leq b$ for all $c\in C$ but $x\cdot
y>b$.

\meno
{\bf Proof.}\quad
Let $d$ be the greatest lower bound of $\Vert x-c\Vert^2$ for all
$c\in C$. Then there's a sequence of vectors~$c_k$ with $\Vert
x-c_k\Vert^2 <d+1/k$; this infinite set of vectors contained in the
sphere $\{\,y\mid\Vert x-y\Vert^2\leq d+1\,\}$ 
must have a limit point~$c_{\infty}$,
and $c_{\infty}\in C$ since $C$ is closed. Therefore $\Vert
x-c_{\infty}\Vert^2\geq d$; in fact $\Vert x-c_{\infty}\Vert^2=d$,
since $\Vert x-c_{\infty}\Vert\leq\Vert x-c_k\Vert+\Vert
c_k-c_{\infty}\Vert$ and the right-hand side can be made arbitrarily
close to~$d$. Since $x\notin C$, we must have $d>0$. Now let
$y=x-c_{\infty}$ and $b=c_{\infty}\cdot y$. Clearly $x\cdot y=y\cdot
y+b>b$. And if $c$ is any element of~$C$ and $\epsilon$ is any small
positive number, the vector $\epsilon c+(1-\epsilon)c_{\infty}$ is
in~$C$, hence 
$\bigl\Vert x-\bigl(\epsilon
c+(1-\epsilon)c_{\infty}\bigr)\bigr\Vert^2\geq d$. But
$$\eqalign{\bigl\Vert x-\bigl(\epsilon
c+(1-\epsilon)c_{\infty}\bigr)\bigr\Vert^2-d
&=\Vert x-c_{\infty}-\epsilon(c-c_{\infty})\Vert^2-d\cr
\noalign{\smallskip}
&=\null-2\epsilon y\cdot(c-c_{\infty})+\epsilon^2\,\Vert
c-c_{\infty}\Vert^2\cr}$$ 
can be nonnegative for all small $\epsilon$ only if $y\cdot
(c-c_{\infty})\leq 0$, i.e., $c\cdot y\leq b$.\ \pfbox

\meno
If $A$ is any set of vectors, let $A^{\ast}=\{\,b\mid a\cdot b\geq
0\hbox{ for all $a\in A$}\,\}$. 

The following facts are immediate:

\meno
{\bf F2.}\quad
$A\subseteq A'$ implies $A^{\ast}\supseteq A'^{\ast}$.

\meno
{\bf F3.}\quad
$A\subseteq A^{\ast\ast}$.

\meno
{\bf F4.}\quad
$A^{\ast}$ is a closed cone.

\medskip
From F1 we also get a result which, in the special case that $C=\{\,Ax\mid
x\ge0\,\}$ for a matrix~$A$, is called Farkas's Lemma:

\meno
{\bf F5.}\quad
If $C$ is a closed cone, $C=C^{\ast\ast}$.

\meno
{\bf Proof.}\quad
Suppose $x\in C^{\ast\ast}$ and $x\notin C$, and let $(y,b)$ be a
separating hyperplane as in~F1. Then $(y,0)$ is also a
separating hyperplane; for we have $x\cdot y>b\geq 0$ because $0\in
C$, and we cannot have $c\cdot y>0$ for $c\in C$ because $(\lambda
c)\cdot y$ would then be unbounded. But then $c\cdot(-y)\geq 0$ for
all $c\in C$, so $-y\in C^{\ast}$; hence $x\cdot(-y)\geq 0$,
a~contradiction. \ \pfbox

\medskip
If $A$ and $B$ are sets of vectors, we define
 $A+B=\{\,a+b\mid a\in A\ {\rm and}\ b\in B\,\}$.

\meno
{\bf F6.}\quad
If $C$ and $C'$ are closed cones, $(C\cap
C')^{\ast}=C^{\ast}+C'^{\ast}$.

\meno
{\bf Proof.}\quad
If $A$ and $B$ are arbitrary sets we have $A^{\ast}+B^{\ast}\subseteq
(A\cap B)^{\ast}$, for if $x\in A^{\ast}+B^{\ast}$ and $y\in A\cap B$
then $x\cdot y=a\cdot y+b\cdot y\geq 0$. If $A$ and~$B$ are arbitrary
sets including~0 then $(A+B)^{\ast}\subseteq A^{\ast}\cap B^{\ast}$ 
by~F2, because $A+B\supseteq A$ and $A+B\supseteq B$. Thus for
arbitrary~$A$ and~$B$ we have $(A^{\ast}+B^{\ast})^{\ast}\subseteq
A^{\ast\ast}\cap B^{\ast\ast}$, hence
$$(A^{\ast}+B^{\ast})^{\ast\ast}\supseteq (A^{\ast\ast}\cap
B^{\ast\ast})^{\ast}\,.$$ 
Now let $A$ and $B$ be closed cones; apply F5 to get
$A^{\ast}+B^{\ast}\supseteq(A\cap B)^{\ast}$.\ \pfbox

\meno
{\bf F7.}\quad
If $C$ and $C'$ are closed cones, $(C+C')^{\ast}=C^{\ast}\cap
C'^{\ast}$. (I~don't need this but I~might as well state it.) 
{\bf Proof.} F6~says $(C^{\ast}\cap
C'^{\ast})^{\ast}=C^{\ast\ast}+C'^{\ast\ast}$; apply F5 and $\ast$
again. \ \pfbox

\meno
{\bf F8.}\quad
Let $S$ be any set of indices and let $A_S=\{\,a\mid a_s=0\hbox{ for
all }s\in S\,\}$, and let $\overline{S}$ be all the indices not in~$S$.
 Then
$$A_S^{\ast}=A_{\overline{S}}\,.$$

\meno
{\bf Proof.}\quad
If $b_s=0$ for all $s\notin S$ and $a_s=0$ for all $s\in S$, obviously
$a\cdot b=0@$; so $A_{\overline{S}}\subseteq A_S^{\ast}$. If $b_s\neq 0$ for
some $s\notin S$ and $a_t=0$ for all $t\neq s$ and $a_s=-b_s$ then
$a\in A_S$ and $a\cdot b <0@$; so $b\notin A_S^{\ast}$, hence
$A_{\overline{S}}\supseteq A_S^{\ast}$.\ \pfbox

\meno
{\bf 9. Definite proof of a semidefinite fact.}\quad
Now we are almost ready to prove the result needed in the proof of
Lemma~7.

Let $D$ be the set of real symmetric positive semidefinite matrices
(called ``spuds'' henceforth for brevity),
 considered as vectors in $N$-dimensional space where
$N={1\over 2}(n+1)n$. We use the inner product $A\cdot B={\rm
tr}\,A^TB$; this is justified if we divide off-diagonal elements
by~$\sqrt{2}$. For example, if $n=3$ the correspondence between
6-dimensional vectors and $3\times 3$ symmetric matrices~is
$$(a,b,c,d,e,f)\ \leftrightarrow\ \pmatrix{a&d/\sqrt{2}&e/\sqrt{2}\cr
\noalign{\smallskip}
d/\sqrt{2}&b&f/\sqrt{2}\cr
\noalign{\smallskip}
e/\sqrt{2}&f/\sqrt{2}&c\cr}$$
preserving sum, scalar product, and dot product. Clearly $D$ is a
closed cone.

\meno
{\bf F9.}\quad
$D^{\ast}=D$.

\meno
{\bf Proof.}\quad
If $A$ and $B$ are spuds then $A=X^TX$ and $B=Y^TY$ and $A\cdot B={\rm
tr}\,X^TX\,Y^TY={\rm tr}\,XY^TYX^T=(YX^T)\cdot(YX^T)\geq 0@$; hence
$D\subseteq D^{\ast}$. (In fact, this argument shows that $A\cdot B=0$ iff
$AB=0$, for any spuds $A$ and~$B$, since $A=A^T$.)

If $A$ is symmetric but has a negative eigenvalue $\lambda$ we can
write 
$$A=Q\,{\rm diag}\,(\lambda,\lambda_2,\ldots,\lambda_n)\,
%\pmatrix{\lambda\cr &\lambda_2\cr &&\ddddots\cr &&&\lambda_n\cr}
Q^T$$
 for some orthogonal matrix~$Q$. Let
$B=Q\,{\rm diag}\,(1,0,\ldots,0)Q^T$; then $B$ is a spud,
and
$$A\cdot B={\rm tr}\;A^TB={\rm tr}\;Q\,{\rm diag}\,(\lambda,0,\ldots,0)\,
Q^T=\lambda <0\,.$$
So $A$ is not in $D^{\ast}$; this proves $D\supseteq D^{\ast}$.\
\pfbox

\medskip
Let $E$ be the set of all real symmetric matrices such that $E_{uv}=0$
when $u\adj v$ in a graph~$G$;  
let $F$ be the set of all real symmetric matrices such that $F_{uv}=0$
when $u=v$ or $u\noadj v$.
The Fact  stated in Section~7 is now equivalent in our new notation to

\meno
{\bf Fact.}\qquad$(D\cap E)^{\ast}\subseteq D+F$.

But we know that 
$$\vcenter{\halign{$#$\hfil\ &$#$\hfil\qquad&#\hfil\cr
(D\cap E)^{\ast}&=D^{\ast}+E^{\ast}&by F6\cr
\noalign{\smallskip}
&=D+F&by F9 and F8.\ \pfbox\cr}}$$

\meno
{\bf 10. Another characterization.}\quad
Remember $\vartheta$, $\vartheta_1$, $\vartheta_2$, and $\vartheta_3$?
We are now going to introduce yet another function
$$\vartheta_4(G,w)=\max\left\{\,\sum_v c(b_v)w_v\biggm| b\hbox{ is an
orthogonal labeling of }\overline{G}\,\right\}\,.\eqno(10.1)$$

\proclaim
Lemma. $\vartheta_3(G,w)\leq \vartheta_4(G,w)$.

\noindent
{\bf Proof.}\quad
Suppose $b$ is a normalized orthogonal labeling of~$\overline{G}$ that
achieves the maximum~$\vartheta_3$; and suppose the vectors of this
labeling have dimension~$d$. Let
$$x_k=\sum_vb_{kv}\sqrt{w_v}\,,\qquad{\rm for}\ 1\leq k\leq
d\,;\eqno(10.2)$$ 
then
$$\vartheta_3(G,w)=\sum_{u,v}\sqrt{w_u}\,b_u\cdot b_v\,\sqrt{w_v}
=\sum_{u,v,k}\sqrt{w_uw_v}\,b_{ku}b_{kv}=\sum_kx_k^2\,.$$
Let $Q$ be an orthogonal $d\times d$ matrix whose first row is
$(x_1/\sqrt{\vartheta_3},\ldots,x_d/\sqrt{\vartheta_3})^T$, and let
$b'_v=Qb_v$. Then $b'_u\cdot b'_v=b_u^TQ^TQb_v=b_u^Tb_v=b_u\cdot b_v$,
so $b'$ is a normalized orthogonal labeling of~$\overline{G}$. Also
$$\eqalignno{x'_k=\sum_vb'_{kv}\sqrt{w_v}
&=\sum_{v,j}Q_{kj}b_{jv}\sqrt{w_v}\cr
\noalign{\smallskip}
&=\sum_jQ_{kj}x_j=\cases{\sqrt{\vartheta_3}\,,&$k=1$;\cr
\noalign{\smallskip}
0\,,&$k>1$.\cr}&(10.3)\cr}$$
Hence by Cauchy's inequality
$$\eqalignno{\vartheta_3(G,w)=\biggl(\sum_vb'_{1v}\sqrt{w_v}\biggr)^2
&\leq\biggl(\sum_v\Vert b'_v\Vert^2\biggr)\biggl(
\sum_{\scriptstyle v\atop \scriptstyle b'_v\neq
0}\;{b'^2_{1v}\over\Vert b'_v\Vert^2}\,w_v\biggr)\cr
\noalign{\smallskip}
&=\sum_v c(b'_v)w_v\leq\vartheta_4(G,w)&(10.4)\cr}$$
because $\sum_v\Vert b'_v\Vert^2=\sum_v\Vert b_v\Vert^2=1$. \ \pfbox

\meno
{\bf 11. The final link.}\quad
Now we can close the loop:

\proclaim
Lemma. $\vartheta_4(G,w)\leq \vartheta(G,w)$.

\noindent
{\bf Proof.}\quad
If $b$ is an orthogonal labeling of $\overline{G}$ that achieves the
maximum~$\vartheta_4$, we will show that the real labeling~$x$
defined by $x_v=c(b_v)$ is in ${\tt TH}(G)$. Therefore
$\vartheta_4(G,w)=w\cdot x$ is $\leq \vartheta(G,w)$.

We will prove that if $a$ is any orthogonal labeling of $G$, and if
$b$ is any orthogonal labeling of~$\overline{G}$, then
$$\sum_v c(a_v)@c(b_v)\leq 1\,.\eqno(11.1)$$
Suppose $a$ is a labeling of dimension $d$ and $b$ is of
dimension~$d'$. Then consider the $d\times d'$ matrices
$$A_v=a_v\0b_v^T\eqno(11.2)$$
as elements of a vector space of dimension~$dd'$. If $u\neq v$ we have
$$A_u\cdot A_v={\rm tr}\;A_u^TA_v\0={\rm tr}\;b_u\0a_u^Ta_v\0b_v^T={\rm
tr}\; a_u^Ta_v\0b_v^Tb_u\0=0\,,\eqno(11.3)$$
because $a_u^Ta\0_v=0$ when $u\noadj v$ and $b_v^Tb\0_u=0$ when $u\adj v$.
If $u=v$ we have
$$A_v\cdot A_v=\Vert a_v\Vert^2\,\Vert b_v\Vert^2\,.$$
The upper left corner element of $A_v$ is $a_{1v}b_{1v}$, hence the
``cost'' of~$A_v$ is $(a_{1v}b_{1v})^2/\,\Vert
A_v\Vert^2=c(a_v)@c(b_v)$. This, with (11.3), proves (11.1). (See the 
proof of Lemma~1.)\ \pfbox

\meno
{\bf 12. The main theorem.}\quad
Lemmas 5, 6, 7, 10, and 11 establish the five inequalities claimed in
(5.1); hence all five variants of~$\vartheta$ are the same function
of~$G$ and~$w$. Moreover, all the inequalities in those five proofs
are equalities $\bigl($with the exception of (11.1)$\bigr)$. We can
summarize the results as follows.

\proclaim
Theorem. For all graphs $G$ and any nonnegative real labeling $w$
of~$G$ we have
$$\vartheta(G,w)=\vartheta_1(G,w)=
\vartheta_2(G,w)=\vartheta_3(G,w)=\vartheta_4(G,w)\,.\eqno(12.1)$$
Moreover, if $w\neq 0$, there exist orthogonal labelings $a$ and $b$
of~$G$ and $\overline{G}$, respectively, such that
$$\eqalignno{c(a_v)&=w_v/\vartheta\,;&(12.2)\cr
\noalign{\smallskip}
\sum c(a_v)@c(b_v)&=1\,.&(12.3)\cr}$$

\noindent
{\bf Proof.}\quad
Relation (12.1) is, of course, (5.1); and (12.2) is (5.3). The desired
labeling~$b$ is what we called~$b'$ in the proof of Lemma~10. The
fact that the application of Cauchy's inequality in (10.4) is
actually an equality,
$$\vartheta=\biggl(\sum_vb_{1v}\sqrt{w_v}\,\biggr)^2
=\biggl(\sum_v\,\Vert b_v\Vert^2\biggr)\biggl(\,
\sum_{\scriptstyle v\atop\scriptstyle b_v\neq
0}\; {b^2_{1v}\over\Vert b_v\Vert^2}\, w_v\biggr)\,,\eqno(12.4)$$
tells us that the vectors whose dot product has been squared are
proportional: There is a number~$t$ such that
$$\Vert b_v\Vert =t\,{b_{1v}\sqrt{w_v}\over\Vert b_v\Vert}\,,\quad
{\rm if}\ b_v\neq 0\,;\qquad\Vert b_v\Vert =0\quad{\rm iff} \quad
b_{1v}\sqrt{w_v}=0\,.\eqno(12.5)$$
The labeling in the proof of Lemma 10 also satisfies
$$\sum_v\,\Vert b_v\Vert^2=1\,;\eqno(12.6)$$
hence $t=\pm 1/\sqrt{\vartheta}\,$. 

We can now show 
$$c(b_v)=\Vert b_v\Vert^2\,\vartheta/w_v\,,\quad{\rm when}\ w_v\neq
0\,.\eqno(12.7)$$
This relation is obvious if $\Vert b_v\Vert=0@$; otherwise we have
$$c(b_v)={b^2_{1v}\over\Vert b_v\Vert^2}={\Vert b_v\Vert^2\over
t^2w_v}$$
by (12.5). Summing the product of 
(12.2) and (12.7) over $v$ gives (12.3).\ \pfbox

\meno
{\bf 13. The main converse.}\quad
The nice thing about Theorem 12 is that conditions (12.2) and (12.3)
also provide a {\it certificate\/} that a given value~$\vartheta$ is
the minimum or maximum stated in the definitions of $\vartheta$,
$\vartheta_1$, $\vartheta_2$, $\vartheta_3$, and $\vartheta_4$.

\proclaim
Theorem. If $a$ is an orthogonal labeling of $G$ and $b$ is an
orthogonal labeling of\/ $\overline{G}$ such that relations (12.2) and
(12.3) hold for some~$\vartheta$ and~$w$, then $\vartheta$ is the
value of $\vartheta(G,w)$.

\noindent
{\bf Proof.}\quad
Plugging (12.2) into (12.3) gives $\sum w_vc(b_v)=\vartheta$, hence
$\vartheta\leq\vartheta_4(G,w)$ by definition of~$\vartheta_4$. Also,
$$\max_v\;{w_v\over c(a_v)}=\vartheta\,,$$
hence $\vartheta\geq \vartheta_1(G,w)$ by definition of~$\vartheta_1$.
\ \pfbox

\meno
{\bf 14. Another look at {\tt TH}.}\quad
We originally defined $\vartheta(G,w)$ in (4.1) in terms of the convex
set~{\tt TH} defined in section~2:
$$\vartheta(G,w)=\max\{\,w\cdot x\mid x\in {\tt
TH}(G)\,\}\,,\qquad{\rm when}\ w\geq 0\,.\eqno(14.1)$$
We can also go the other way, defining {\tt TH} in terms
of~$\vartheta$:
$${\tt TH}(G)=\{\,x\geq 0\mid w\cdot x\leq\vartheta(G,w)\hbox{ for all
}w\geq 0\,\}\,.\eqno(14.2)$$
Every $x\in {\tt TH}(G)$ belongs to the right-hand set, by (14.1).
Conversely, if $x$ belongs to the right-hand set and if $a$ is any
orthogonal labeling of~$G$, not entirely zero, let $w_v=c(a_v)$, so
that $w\cdot x=\sum_v c(a_v)x_v$. Then
$$\vartheta_1(G,w)\leq\max_v\bigl(w_v/c(a_v)\bigr)=1$$
by definition (5.2), so we know by Lemma 5 that $\sum c(a_v)x_v\leq
1$. This proves that $x$ belongs to ${\tt TH}(G)$.

Theorem 12 tells us even more.

\proclaim
Lemma. ${\tt TH}(G)=\{\,x\geq 0\mid\vartheta(\overline{G},x)\leq 1\,\}$.

\noindent
{\bf Proof.}\quad
By definition (10.1),
$$\vartheta_4(\overline{G},w)=\max\left\{\,\sum_v c(a_v)w_v\mid a\hbox{ is
an orthogonal labeling of }G\,\right\}\,.\eqno(14.3)$$
Thus $x\in {\tt TH}(G)$ iff $\vartheta_4(\overline{G},x)\leq 1$, when
$x\geq 0$.\ \pfbox

\proclaim
Theorem. ${\tt TH}(G)=\{\,x\mid x_v=c(b_v)\hbox{ for some orthogonal
labeling $b$ of $\overline{G}$}\,\}$.

\noindent
{\bf Proof.}\quad
We already proved in (11.1) that the right side is contained in the
left. 

Let $x\in {\tt TH}(G)$ and let $\vartheta=\vartheta(\overline{G},x)$. By
the lemma, $\vartheta\leq 1$. 
Therefore, by (12.2), there is an orthogonal labeling~$b$ of~$\overline{G}$
such that $c(b_v)=x_v/\vartheta\geq x_v$ for all~$v$. It is easy to reduce the
cost of any vector in an orthogonal labeling to any desired value, simply by
increasing the dimension
 and giving this vector an appropriate nonzero value in
the new component while all other vectors remain zero there. The dot products
are unchanged, so the new labeling is still orthogonal. Repeating this
construction for each~$v$ produces a labeling with $c(b_v)=x_v$. \ \pfbox

\medskip
This theorem makes the definition of $\vartheta_4$ in (10.1) identical to the
definition of~$\vartheta$ in~(4.1).

\meno
{\bf 15. Zero weights.}\quad
Our next result shows that when a weight is zero, the corresponding vertex
might as well be absent from the graph. 

\proclaim
Lemma. Let $U$ be a subset of the vertices~$V$ of a graph~$G$, and let
$G'=G\vert U$ be the graph induced by~$U$ (i.e., the graph on
vertices~$U$ with $u\adj v$ in~$G'$ iff $u\adj v$ in~$G$). Then if $w$ and
$w'$ are nonnegative labelings of~$G$ and~$G'$ such that
$$w_v=w'_v\qquad\hbox{when}\quad v\in U\,,\qquad
w_v=0\qquad\hbox{when}\quad v\notin U\,,\eqno(15.1)$$
we have
$$\vartheta(G,w)=\vartheta(G',w')\,.\eqno(15.2)$$

\noindent
{\bf Proof.}\quad
Let $a$ and $b$ satisfy (12.2) and (12.3) for $G$ and~$w$. 
Then $c(a_v)=0$ for $v\notin U$, so $a\vert U$ and 
$b\vert U$ satisfy
(12.2) and (12.3) for~$G'$ and~$w'$. (Here
$a\vert U$ means the vectors~$a_v$ for $v\in U$.) By Theorem~13, they
determine the same~$\vartheta$.\ \pfbox

\meno
{\bf 16. Nonzero weights.}\quad
We can also get some insight into the significance of nonzero weights
by ``splitting'' vertices instead of removing them.

\proclaim
Lemma. Let $v$ be a vertex of $G$ and let $G'$ be a graph obtained
from~$G$ by adding a new vertex~$v'$ and new edges
$$u\adj v'\quad{\rm iff}\quad u\adj v\,.\eqno(16.1)$$
Let $w$ and $w'$ be nonnegative labelings of~$G$ and~$G'$ such that
$$\eqalignno{w_u&=w'_u\,,\qquad\hbox{when }u\neq v\,;&(16.2)\cr
\noalign{\smallskip}
w_v&=w'_v+w'_{v'}\,.&(16.3)\cr}$$
Then
$$\vartheta(G,w)=\vartheta(G',w')\,.\eqno(16.4)$$

\noindent
{\bf Proof.}\quad
By Theorem 12 there are labelings $a$ and $b$ of~$G$ and~$\overline{G}$
satisfying (12.2) and (12.3). We can modify them to obtain labelings
$a'$ and~$b'$ of~$G'$ and~$\overline{G'}$ as follows, with the vectors
of~$a'$ having one more component than the vectors of~$a$:
$$a'_u={a_u\choose 0}\,,\qquad
b'_u=b_u\,,\qquad\hbox{when }u\neq v\,;\eqno(16.5)$$
\vskip-10pt
$$a'_v={a_v\choose\alpha}\,,\qquad a'_{v'}={a_v\choose -\beta}\,,\qquad
\alpha=\sqrt{{w'_{v'}\over w'_v}}\; \Vert a_v\Vert\,,\qquad
\beta=\sqrt{{w'_v\over w'_{v'}}}\;\Vert a_v\Vert\,;\eqno(16.6)$$
\vskip-10pt
$$b'_v=b'_{v'}=b_v\,.\eqno(16.7)$$
(We can assume by Lemma 15 that $w'_v$ and $w'_{v'}$ are nonzero.) All
orthogonality relations are preserved; and since $v\noadj v'$ in~$G'$,
we also need to verify
$$a'_v\cdot a'_{v'}=\Vert a_v\Vert^2-\alpha\beta =0\,.$$
We have
$$c(a'_v)={c(a_v)\;\Vert a_v\Vert^2\over \Vert a_v\Vert^2+\alpha^2}
={c(a_v)\over 1+w'_{v'}/w'_v}={c(a_v)w'_v\over
w_v}={w'_v\over\vartheta}\,,$$
and similarly $c(a'_{v'})=w'_{v'}/\vartheta$; thus (12.2) and (12.3)
are satisfied by~$a'$ and~$b'$ for~$G'$ and~$w'$. \ \pfbox

\medskip
Notice that if all the weights are integers we can apply this lemma
repeatedly to establish that
$$\vartheta(G,w)=\vartheta(G')\,,\eqno(16.8)$$
where $G'$ is obtained from~$G$ by replacing each vertex~$v$ by a
cluster of~$w_v$ mutually nonadjacent vertices that are adjacent to
each of~$v$'s neighbors. $\bigl($Recall that
$\vartheta(G')=\vartheta(G',\ones)$, by definition (4.2).$\bigr)$
In particular, if $G$ is the trivial graph~$K_2$ and if we assign the
weights $M$ and~$N$, we have $\vartheta\bigl(K_2,(M, 
N)^T\bigr)=\vartheta(K_{M,N})$ where $K_{M,N}$ denotes the complete
bipartite graph on~$M$ and~$N$ vertices.

A similar operation called ``duplicating'' a vertex has a similarly simple
effect:

\proclaim
Corollary. Let $G'$ be constructed from $G$ as in the lemma but with an
additional edge between~$v$ and~$v'$. Then $\vartheta(G,w)=\vartheta(G',w')$ 
if $w'$ is defined by (16.2) and
$$w_v=\max(w'_v,w'_{v'})\,.\eqno(16.9)$$

\noindent
{\bf Proof.}\quad
We may assume that $w_v=w'_v$ and $w'_{v'}\neq 0$. Most of the construction
(16.5)--(16.7) can be used again, but we set $\alpha=0$ and $b'_{v'}=0$ and
$$\beta=\sqrt{{w_v-w'_{v'}\over w'_{v'}}}\;\Vert a_v\Vert\,.$$
Once again the
 necessary and sufficient conditions are readily verified. \ \pfbox

\medskip
If the corollary is applied repeatedly, it tells us that $\vartheta(G)$ is
unchanged when we replace the vertices of~$G$ by cliques.

\meno
{\bf 17. Simple examples.}\quad
We observed in section 4 that $\vartheta(G,w)$ always is at least
$$\vartheta_{\min}=\vartheta(K_n,w)=\max\{w_1,\ldots,w_n\}\eqno(17.1)$$
and at most
$$\vartheta_{\max}=(\overline{K}_n,w)=w_1+\cdots +w_n\,.\eqno(17.2)$$
What are the corresponding orthogonal labelings?

For $K_n$ the vectors of $a$ have no orthogonal constraints, while the
vectors of~$b$ must satisfy $b_u\cdot b_v=0$ for all $u\neq v$. We can
let $a$ be the two-dimensional labeling
$$a_v={\sqrt{w_v}\choose \sqrt{\vartheta-w_v}}\,,\qquad
\vartheta=\vartheta_{\min}\eqno(17.3)$$
so that $\Vert a_v\Vert^2=\vartheta$ and $c(a_v)=w_v/\vartheta$ as
desired; and $b$ can be one-dimensional,
$$b_v=\cases{(1)\,,&if $v=v_{\max}$\cr
\noalign{\smallskip}
(0)\,,&if $v\neq v_{\max}$\cr}\eqno(17.4)$$
where $v_{\max}$ is any particular vertex that maximizes~$w_v$.
Clearly
$$\sum_v c(a_v)@c(b_v)={c(a_{v_{\max}})\over \vartheta}=
{w_{v_{\max}}\over\vartheta}=1\,.$$

For $\overline{K}_n$ the vectors of $a$ must be mutually orthogonal while
the vectors of~$b$ are unrestricted. We can let the vectors~$a$ be the
columns of any orthogonal matrix whose top row contains the element
$$\sqrt{w_v/\vartheta}\,,\qquad\vartheta=\vartheta_{\max}\eqno(17.5)$$
in column $v$. Then $\Vert a_v\Vert^2=1$ and $c(a_v)=w_v/\vartheta$.
Once again a one-dimensional labeling suffices for~$b$; we can let
$b_v=(1)$ for all~$v$.

\meno
{\bf 18. The direct sum of graphs.}\quad
Let $G=G'+G''$ be the graph on vertices
$$V=V'\cup V''\eqno(18.1)$$
where the vertex sets~$V'$ and $V''$ of~$G'$ and~$G''$ are disjoint,
and where $u\adj v$ in~$G$ if and only if $u,v\in V'$ and $u\adj v$
in~$G'$, or $u,v\in V''$ and $u\adj v$ in~$G''$. In this case
$$\vartheta(G,w)=\vartheta(G',w')+\vartheta(G'',w'')\,,\eqno(18.2)$$
where $w'$ and $w''$ are the sublabelings of~$w$ on vertices of~$V'$
and~$V''$. We can prove (18.2) by constructing orthogonal labelings
$(a,b)$ satisfying (12.2) and (12.3).

Suppose $a'$ is an orthogonal labeling of $G'$ such that
$$\Vert a'_v\Vert^2=\vartheta'\qquad
a'_{1v}=\sqrt{w'_v}\,,\eqno(18.3)$$
and suppose $a''$ is a similar orthogonal labeling of~$G''$. If $a'$
has dimension~$d'$ and~$a''$ has dimension~$d''$, we construct a new
labeling~$a$ of dimension $d=d'+d''$ as follows, where $j'$ runs
from~2 to~$d'$ and $j''$ runs from~2 to~$d''$:
$$\vcenter{\halign{\hfil$#\;$&$#$\hfil\qquad
&\hfil$#\;$&$#$\hfil\cr
&\hbox{ if $v\in V'$}&&\hbox{ if $v\in V''$}\cr
\noalign{\medskip}
a_{1v}&=\sqrt{w_v'}=a'_{1v}\,,&a_{1v}&=\sqrt{w''_v}=a''_{1v}\,,\cr
\noalign{\smallskip}
a_{j'v}&=\sqrt{\vartheta/\vartheta'}\,a'_{j'v}\,,&a_{j'v}&=0\,,\cr
\noalign{\smallskip}
a_{(d'+1)v}&=\sqrt{\vartheta''w'_v/\vartheta'}\,,
&a_{(d'+1)v}&=-\sqrt{\vartheta'w''_v/\vartheta''}\,,\cr
\noalign{\smallskip}
a_{(d'+j'')v}&=0\,,&a_{(d'+j'')v}
&=\sqrt{\vartheta/\vartheta''}\,a''_{j''v}\,.\cr}}\eqno(18.4)$$
Now if $u,v\in V'$ we have
$$a_u\cdot
a_v=\sqrt{w'_uw'_v}+{\vartheta\over\vartheta'}\,\bigl(a'_u\cdot
a'_v-\sqrt{w'_uw'_v}\,\bigr)+{\vartheta''\over\vartheta'}\,\sqrt{w'_uw'_v}=
{\vartheta\over\vartheta'}\,a'_u\cdot a'_v\,;\eqno(18.5)$$
thus $a_u\cdot a_v=0$ when $a'_u\cdot a'_v=0$, and
$$\Vert a_v\Vert^2={\vartheta\over\vartheta'}\,\Vert
a'_v\Vert^2=\vartheta\,.\eqno(18.6)$$ 
It follows that $c(a_v)=w_v/\vartheta$ as desired. A similar
derivation holds for $u,v\in V''$. And if $u\in V'$, $v\in V''$, then
$$a_u\cdot a_v=\sqrt{w'_uw''_v}-\sqrt{w'_uw''_v}=0\,.\eqno(18.7)$$

The orthogonal labeling $b$ of $\overline{G'+G''}$ is much simpler;
we just let $b_v=b'_v$ for $v\in V'$ and $b_v=b''_v$ for $v\in V''$.
Then (12.2) and (12.3) are clearly preserved. This proves (18.2).

There is a close relation between the construction (18.4) and the
construction (16.6), suggesting that we might be able to define
another operation on graphs that generalizes both the splitting and
direct sum operation.

\meno
{\bf 19. The direct cosum of graphs.}\quad
If $G'$ and $G''$ are graphs on disjoint vertex sets~$V'$ and~$V''$ as
in section~18, we can also define
$$G=G'\mathbin{\overline+} G''\ 
\Longleftrightarrow \overline{G}=\overline{G'}+
\overline{G''}\,.\eqno(19.1)$$ 
This means $u\adj v$ in $G$ if and only if either $u\adj v$ in~$G'$ or
$u\adj v$ in~$G''$ or $u$ and~$v$ belong to opposite vertex sets. In
this case
$$\vartheta(G,w)=\max\bigl(\vartheta(G',w'),\;
\vartheta(G'',w'')\bigr)\eqno(19.2)$$
and again there is an easy way to construct $(a,b)$ from $(a',b')$ and
$(a'',b'')$ to prove (19.2). Assume ``without lots of generality''
that
$$\vartheta(G',w')\geq\vartheta(G'',w'')\eqno(19.3)$$
and suppose again that we have (18.3) and its counterpart for~$a''$.
Then we can define
$$\vcenter{\halign{\hfil$#\;$&$#$\hfil\qquad
&\hfil$#\;$&$#$\hfil\cr
&\hbox{ if $v\in V'$}&&\hbox{ if $v\in V''$}\cr
\noalign{\medskip}
a_{1v}&=\sqrt{w_{v'}}=a'_{1v}\,,&a_{1v}&=\sqrt{w_{v''}}=a''_{1v}\,,\cr
\noalign{\smallskip}
a_{j'v}&=a'_{j'v}\,,&a_{j'v}&=0\,,\cr
\noalign{\smallskip}
a_{(d'+1)v}&=0\,,
&a_{(d'+1)v}&=\sqrt{(\vartheta'-\vartheta'')w''_v/\vartheta''}\,,\cr
\noalign{\smallskip}
a_{(d'+j'')v}&=0\,,&a_{(d'+j'')v}&=\sqrt{\vartheta'/\vartheta''}\,a''_{j''v}
\,.\cr}}\eqno(19.4)$$
Now $a_v$ is essentially unchanged when $v\in V'$; and when $u,v\in
V''$ we have
$$a_u\cdot
a_v=\sqrt{w''_uw''_v}+\left(\,{\vartheta'\over\vartheta''}-1\right)\, 
\sqrt{w''_uw''_v}+{\vartheta'\over\vartheta''}\,
\bigl(a''_u\cdot
a''_v-\sqrt{w''_uw''_v}\,\bigr)={\vartheta'\over\vartheta''} \,
a''_u\cdot a''_v\,.\eqno(19.5)$$
Again we retain the necessary orthogonality, and we have
$c(a_v)=w_v/\vartheta$ for all~$v$.

For the $b$'s, we let $b_v=b'_v$ when $v\in V'$ and $b_v=0$ when $v\in
V''$. 

\meno
{\bf 20. A direct product of graphs.}\quad
Now let $G'$ and $G''$ be graphs on vertices~$V'$ and~$V''$ and let
$V$ be the $n=n'n''$ ordered pairs
$$V=V'\times V''\,.\eqno(20.1)$$
We define the `strong product',
$$G=G'\;\ast\;G''\eqno(20.2)$$
on $V$ by the rule
$$\eqalignno{&(u',u'')\adj (v',v'')\quad{\rm or}\quad
(u',u'')=(v',v'')\quad {\rm in} \ G\cr
\noalign{\smallskip}
&\quad\Longleftrightarrow\ (u'\adj v'\ {\rm or}\ u'=v'\ {\rm in}\
G')\quad {\rm and}\quad (u''\adj v''\ {\rm or}\ u''=v''\ {\rm in}\
G'')\,.&(20.3)\cr}$$

In this case we have, for example, $K_{n'}\;\ast\; K_{n''}=K_{n'n''}$
and $\overline{K}_{n'}\;\ast\;\overline{K}_{n''}=\overline{K}_{n'n''}$. More
generally, if
$G'$ is regular of degree $r'$ and $G''$ is regular of degree~$r''$,
then $G'\;\ast\;G''$ is regular of degree
$(r'+1)(r''+1)-1=r'r''+r'+r''$.

I don't know the value of $\vartheta(G,w)$ for arbitrary~$w$, but I do
know it in the special case
$$w_{(v',v'')}=w'_{v'}w''_{v''}\,.\eqno(20.4)$$

\proclaim
Lemma. If $G$ and $w$ are given by (20.2) and (20.4), then
$$\vartheta(G,w)=\vartheta(G',w')\,\vartheta(G'',w'')\,.\eqno(20.5)$$

\noindent
{\bf Proof.}\quad
[\Lov]\enspace
Given orthogonal labelings $(a',b')$ and $(a'',b'')$ of $G'$
and~$G''$, we let $a$ be the Hadamard product
$$a_{(j',j'')(v',v'')}=a'_{j'v'}a''_{j''v''}\,,\qquad
1\leq j'\leq d'\,,\quad 1\leq j''\leq d''\,,\eqno(20.6)$$ 
where $d'$ and $d''$ are the respective dimensions of the 
vectors in $a'$ and~$a''$. Then
$$\eqalignno{a_{(u',u'')}\cdot a_{(v',v'')}
&=\sum_{j',j''}a'_{j'u'}a''_{j''u''}a'_{j'v'}a''_{j''v''}\cr
\noalign{\smallskip}
&=(a'_{u'}\cdot a'_{v'})(a''_{u''}\cdot a''_{v''})\,.&(20.7)\cr}$$
Thus $\Vert a_{(v',v'')}\Vert^2=\Vert a'_{v'}\Vert^2\,\Vert
a''_{v''}\Vert^2$ and 
$$c(a_{(v'v'')})=c(a'_{v'})@c(a''_{v''})\,.\eqno(20.8)$$
The same construction is used for $b$ in terms of $b'$ and $b''$.

All necessary orthogonalities are preserved, because we have
$$\eqalign{%
&(u',u'')\adj (v',v'')\hbox{ and }(u',u'')\neq (v',v'')\hbox{ in }G\cr
\noalign{\vskip3pt}
&\qquad \Rightarrow (u'\adj v'\hbox{ and }u'\neq v'\hbox{ in }G')
\hbox{ or }(u''\adj v''\hbox{ and }u''\neq v''\hbox{ in }G'')\cr
\noalign{\vskip3pt}
&\qquad \Rightarrow b_{(u',u'')}\cdot b_{(v',v'')}=0\,;\cr
\noalign{\smallskip}
&(u',u'')\noadj (v',v'')\hbox{ and }(u',u'')\neq (v',v'')\hbox{ in }G\cr
\noalign{\vskip3pt}
&\qquad \Rightarrow (u'\noadj v'\hbox{ and }u'\neq v'\hbox{ in }G')
\hbox{ or }(u''\noadj v''\hbox{ and }u''\neq v''\hbox{ in
}G'')\cr
\noalign{\vskip3pt}
&\qquad \Rightarrow a_{(u',u'')}\cdot a_{(v',v'')}=0\,.\cr}$$
(In fact one of these relations is $\Leftrightarrow$, but we need only
$\Rightarrow$ to make (20.7) zero when it needs to be zero.) Therefore
$a$ and~$b$ are orthogonal labelings of~$G$ that satisfy (12.2) and
(12.3). \pfbox

\meno
{\bf 21. A direct coproduct of graphs.}\quad
Guess what? We also define
$$G=G'\;\overline{\ast}\;G''\ \Longleftrightarrow\ \overline{G}=
\overline{G'}\;\ast\;\overline{G''}\,.\eqno(21.1)$$ 
This graph tends to be ``richer'' than $G'\;\ast\;G''$; we have
$$\eqalignno{&(u',u'')\adj (v',v'')\hbox{ and }(u',u'')\neq (v',v'')
\hbox{ in }G\cr
&\qquad \Longleftrightarrow (u'\adj v'\hbox{ and }u'\neq v'
\hbox{ in }G')\hbox{ or }(u''\adj v''\hbox{ and }u''\neq v''
\hbox{ in }G'')\,.&(21.2)\cr}$$
Now, for instance, if $G'$ is regular of degree~$r'$ and $G''$ is
regular of degree~$r''$, then
$$G'\;\bar{\ast}\;G''\hbox{ is regular of degree }
n'n''-(n'-r')(n''-r'')=r'n''+r''n'-r'r''\,.$$
(This is always $\geq r'r''+r'+r''$, because
$r'(n''-1-r'')+r''(n'-1-r')\geq 0$.)
Indeed, $G'\;\overline{\ast}\;G''\supseteq G'\;\ast\;G''$ for all graphs~$G'$
and~$G''$. 
The Hadamard product construction used in section~20 can be applied
word-for-word to prove that
$$\vartheta(G,w)=\vartheta(G',w')\,\vartheta(G'',w'')\eqno(21.3)$$
when $G$ satisfies (21.1) and $w$ has the special factored form
(20.4).

It follows that many graphs have identical $\vartheta$'s:

\proclaim
Corollary. If $G'\,\ast\,G''\subseteq G\subseteq G'\;\overline{\ast}\; G''$
and $w$ satisfies (20.4), then (21.3) holds.

\noindent
{\bf Proof.}\quad
This is just the monotonicity relation (4.3). The reason it works is that we
have $a_{(u',v')}\cdot a_{(u'',v'')}=b_{(u',v')}\cdot b_{(u'',v'')}=0$ for 
all pairs of vertices $(u',u'')$ and $(v',v'')$ whose adjacency differs in
$G'\,\ast\,G''$ and $G'\;\overline{\ast}\;G''$. \ \pfbox

\medskip
Some small examples will help clarify the results of the past few sections. 
Let $P_3$ be the path of length~2 on 3~vertices,
$\bullet\!$---$\!\bullet\!$---$\!\bullet$, 
and consider the four graphs we get by
taking its strong product and coproduct with~$\overline{K}_2$ and~$K_2$:
$$\leftline{\qquad$\overline K_2\,\ast\,P_3 = \mpfig1\qquad\qquad
  \vartheta=\max(u+w,v)+\max(x+z,y)$}$$
$\bigl($Since $P_3$ may be regarded as $\overline{K}_2
\mathbin{\overline +} K_1$ and
$\overline{K}_2$ is $K_1+K_2$, this graph is 
$$\bigl((K_1+K_1)
\mathbin{\overline +} K_1\bigr)+\bigl((K_1+K_1)
\mathbin{\overline +}K_1\bigr)$$
and the formula for $\vartheta$ follows from (18.2) and (19.2).$\bigr)$
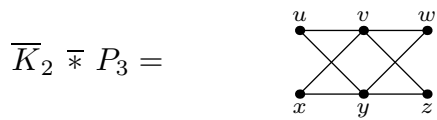
$$\leftline{\qquad$\overline K_2\,\mathbin{\overline\ast}\,P_3 = \mpfig2
\qquad\qquad  \vartheta=\max(u+w+x+z,v+y)$}$$
(This graph is $\overline{K}_2\mathbin{\overline +}\overline{K}_4$; 
we could also obtain
it by applying Lemma~16 three times to~$P_3$.)
$$\leftline{\qquad$K_2\,\ast\,P_3 = \mpfig3\qquad\qquad
  \vartheta=\max\bigl(\max(u,x)+\max(w,z),\max(v,y)\bigr)$}
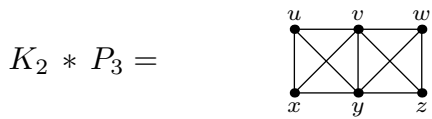
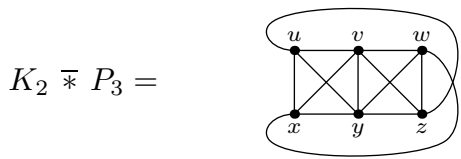
$$
$$\leftline{\qquad$K_2\,\mathbin{\overline\ast}\,P_3 = \mpfig4\qquad\qquad
  \vartheta=\max\bigl(\max(u+w,x+z),\max(v,y)\bigr)$}$$
If the weights satisfy $u=\lambda x$, $v=\lambda y$, $w=\lambda z$ for some
parameter~$\lambda$, the first two formulas for~$\vartheta$ both reduce to
$(1+\lambda)\max(u+w,v)$, in agreement with (20.5) and (21.3). Similarly, the
last two formulas for~$\vartheta$ reduce to $\max(1,\lambda)\max(u+w,v)$ in
such a case.

\meno
{\bf 22. Odd cycles.}\quad
Now let $G=C_n$ be the graph with vertices $0,1,\ldots,n-1$ and
$$u\adj v\ \Longleftrightarrow u-v\equiv \pm
1\pmod{n}\,,\eqno(22.1)$$
where $n$ is an {\it odd\/} number. A~general formula for
$\vartheta(C_n,w)$ appears to be very difficult; but we can compute
$\vartheta(C_n)$ without too much labor when all weights are~1, because of
the cyclic symmetry.

It is easier to construct orthogonal labelings of $\overline{C}_n$
than of~$C_n$, so we begin with that. Given a vertex~$v$, $0\leq v<n$,
let $b_v$ be the three-dimensional vector
$$b_v=\pmatrix{\alpha\cr
\cos v\varphi\cr
\sin v\varphi\cr},\eqno(22.2)$$
where $\alpha$ and $\varphi$ remain to be determined. We have
$$\eqalignno{b_u\cdot b_v&=\alpha^2+\cos u\varphi\;\cos v\varphi+\sin
u\varphi\;\sin v\varphi\cr
&=\alpha^2+\cos (u-v)\varphi\,.&(22.3)\cr}$$
Therefore we can make $b_u\cdot b_v=0$ when $u\equiv v\pm 1$ by
setting
$$\alpha^2=-\cos\varphi\,,\qquad\varphi ={\pi(n-1)\over
n}\,.\eqno(22.4)$$ 
This choice of $\varphi$ makes $n\varphi$ a multiple of~$2\pi$,
because $n$ is odd. We have found an orthogonal labeling~$b$
of~$\overline{C}_n$ such that
$$c(b_v)={\alpha^2\over 1+\alpha^2}={\cos\pi/n\over 1+\cos
\pi/n}\,.\eqno(22.5)$$ 

Turning now to orthogonal labelings of $C_n$, we can use $(2n-1)$-dimensional
vectors
$$a_v=\pmatrix{\alpha_0\cr
\noalign{\vskip3pt}
\alpha_1\cos v\varphi\cr
\noalign{\vskip3pt}
\alpha_1\sin v\varphi\cr
\noalign{\vskip3pt}
\alpha_2\cos 2v\varphi\cr
\noalign{\vskip3pt}
\alpha_2\sin 2v\varphi\cr
\vdots\cr
\noalign{\vskip3pt}
\alpha_{n-1}\cos (n-1)v\varphi\cr
\noalign{\vskip3pt}
\alpha_{n-1}\sin (n-1)v\varphi\cr}\,,\eqno(22.6)$$
with $\varphi=\pi(n-1)/n$ as before. As in (22.3), we find
$$a_u\cdot a_v=\sum_{k=0}^{n-1}\alpha_k^2\cos k(u-v)\varphi\,;\eqno(22.7)$$
so the result depends only on $(u-v)\bmod n$. 
Let $\omega=e^{i\varphi}$. We can
find values of~$\alpha_k$ such that 
$a_u\cdot a_v=x_{(u-v)\bmod n}$ by solving the equations
$$x_j=\sum_{k=0}^{n-1}\alpha_k^2\omega^{jk}\,.\eqno(22.8)$$
Now $\omega$ is a primitive $n$th root of unity; i.e., $\omega^k=1$ iff $k$ is
a multiple of~$n$. So (22.9) is just a finite Fourier transform, and we can
easily invert it: For $0\leq m<n$ we have
$$\sum_{j=0}^{n-1}\omega^{-mj}x_j
=\sum_{k=0}^{n-1}\alpha_k^2\sum_{j=0}^{n-1}\omega^{j(k-m)}
=n\,\alpha_m^2\,.$$
In our case we want a solution with $x_2=x_3=\cdots =x_{n-2}=0$,
 and we can set
$x_0=1$, $x_{n-1}=x_1=x$, so we find
$$n\,\alpha_k^2=x_0+\omega^{-k}x_1+\omega^kx_{n-1}=1+2x\cos k\varphi\,.$$
We must choose $x$ so that these values are nonnegative; this means $2x\leq
-1/\cos\varphi$, since $\cos k\varphi$ is most negative when $k=1$. 
Setting $x$ to this maximum value yields
$$c(a_v)=\alpha_0^2={1\over n}\;\left(1-{1\over\cos\varphi}\right)=
{1+\cos\pi/n\over n\cos\pi/n}\,.\eqno(22.9)$$
So (22.5) and (22.9) give
$$\sum_vc(a_v)@c(b_v)=\sum_v 1/n=1\,.\eqno(22.10)$$
This is (12.3), hence from (12.2) we know that
$\vartheta(C_n)=\lambda$. We have proved, in fact, that
$$\eqalignno{\vartheta(C_n,\ones)&={n\cos \pi/n\over 1+\cos\pi/n}\,;
&(22.11)\cr
\noalign{\smallskip}
\vartheta(\overline{C}_n,\ones)&={1+\cos\pi/n\over\cos
\pi/n}\,.&(22.12)\cr}$$
When $n=3$, $C_n=K_n$ and these values agree with $\vartheta(K_3)=1$,
$\vartheta(\overline{K}_3)=3$; when $n=5$, $\overline{C}_5$~is 
isomorphic to~$C_5$ so
$\vartheta(C_5)=\sqrt{5}$; when $n$ is large,
$$\vartheta(C_n)={n\over 2}-{\pi^2\over 8n}+O(n^{-3})\,;
\qquad
\vartheta(\overline{C}_n)=2+{\pi^2\over 2n^2}+O(n^{-4})\,.\eqno(22.13)$$

Instead of an explicit construction of vectors $a_v$ as in (22.6), we could
also find $\vartheta(C_n)$ by using the matrix characterization~$\vartheta_2$ 
of section~6. When all weights are~1, a~feasible $A$ has 1 everywhere 
except on the superdiagonal, the subdiagonal, and the corners. This 
suggests that we look  at ``circulant'' matrices; for example, when $n=5$,
$$A=\pmatrix{1&1+x&1&1&1+x\cr
1+x&1&1+x&1&1\cr
1&1+x&1&1+x&1\cr
1&1&1+x&1&1+x\cr
1+x&1&1&1+x&1\cr}=J+xP+xP^{-1}\,,\eqno(22.14)$$
where $J$ is all~1's and $P$ is the permutation matrix taking~$j$ into
$(j+1)\bmod n$. It is well known and not difficult to prove that the
eigenvalues of the circulant matrix $a_0I+a_1P+\cdots +a_{n-1}P^{n-1}$
are
$$\sum_{0\leq j<n}\omega^{kj}a_j\,,\qquad 0\leq k<n\,,\eqno(22.15)$$
where $\omega =e^{2\pi i/n}$.
$\bigl($Indeed, it suffices to find the eigenvalues of~$P$ itself.
This $\omega$ is a different primitive root of unity from the~$\omega$ 
we used in (22.8).$\bigr)$ Hence the eigenvalues of (22.14) are
$$n+2x\,,\quad x(\omega+\omega^{-1})\,,\quad
x(\omega^2+\omega^{-2})\,,\ \ldots\,,\quad
x(\omega^{n-1}+\omega^{1-n})\,.\eqno(22.16)$$

We minimize the maximum of these values if we choose $x$ so that
$$n+2x=-2x\cos \pi/n\,;$$
then
$$\Lambda(A)=-2x\cos \pi/n={n\cos \pi/n\over 1+\cos\pi/n}\eqno(22.17)$$
is the value of $\vartheta(G)$.

If $n$ is even, the graph $C_n$ is bipartite. We will prove later that
bipartite graphs are perfect, hence $\vartheta(C_n)=n/2$ and
$\vartheta(\overline{C}_n)=2$ in the even case.

\meno
{\bf 23. Comments on the previous example.}\quad
The  cycles~$C_n$ provide us with infinitely many graphs~$G$ for
which $\vartheta(G)@\vartheta(\,\overline{G}\,)=n$, and it is natural to
wonder whether this is true in general. Of course it is not: If
$G=\overline{K}_m+K_{n-m}$ then $\overline{G}=K_m\mathbin{\overline +}
 \overline{K_{n-m}}$, hence we
know from Lemmas~18 and~19 that
$$\vartheta(G)=m+1\,,\qquad
\vartheta(\,\overline{G}\,)=\max(1,n-m)\,.\eqno(23.1)$$ 
In particular, we can make
$\vartheta(G)@\vartheta(\,\overline{G}\,)$ as high as 
$n^2\!/4+n/2$ when $m=\lfloor n/2\rfloor$.

We can, however, prove without difficulty that
$\vartheta(G)@\vartheta(\,\overline{G}\,)\geq n$:

\proclaim
Lemma. 
$$\vartheta(G,w)@\vartheta(\,\overline{G},w')\geq w\cdot w'\,.\eqno(23.2)$$

\noindent
{\bf Proof.}\quad
By Theorem 12 there is an orthogonal labeling~$a$ of~$G$ and an
orthogonal labeling~$b$ of~$\overline{G}$ such that
$$c(a_v)=w_v/\vartheta(G,w)\,,\qquad
c(b_v)=w'_v/\vartheta(\,\overline{G},w')\,. \eqno(23.3)$$
By (11.1) we have
$$\sum_v c(a_v)@c(b_v)\leq 1\,.\eqno(23.4)$$
QED.\ \pfbox

\meno
{\bf 24. Regular graphs.}\quad
When each vertex of $G$ has exactly $r$~neighbors, Lov\'asz and
Hoffman observed that the construction in (22.14) can be generalized.
Let $B$ be the adjacency matrix of~$G$, i.e., the $n\times n$ matrix
with
$$B_{uv}=\cases{1\,,&if $u\adj v$;\cr
\noalign{\smallskip}
0\,,&if $u=v$ or $u\noadj v$.\cr}\eqno(24.1)$$

\proclaim
Lemma. If $G$ is a regular graph,
$$\vartheta(G)\leq
{n\Lambda(-B)\over\Lambda(B)+\Lambda(-B)}\,.\eqno(24.2)$$ 

\noindent
{\bf Proof.}\quad
Let $A$ be a matrix analogous to (22.14),
$$A=J+xB\,.\eqno(24.3)$$
Since $G$ is regular, the all-1's vector $\ones$ is an eigenvector
of~$B$, and the other eigenvectors are orthogonal to~$\ones$ so they
are eigenvectors also of~$A$. Thus if the eigenvalues of~$B$ are
$$r=\Lambda(B)=\lambda_1\geq \lambda_2\geq \cdots\geq
\lambda_n=-\Lambda(-B)\,, \eqno(24.4)$$
the eigenvalues of $A$ are
$$n+rx\,,\;x\lambda_2\,,\;\ldots\,,\;x\lambda_n\,.\eqno(24.5)$$
(The Perron-Frobenius theorem tells us that $\lambda_1=r$.) We have
$\lambda_1+\cdots+ \lambda_n={\rm tr}(B)=0$, so $\lambda_n<0$, and we
minimize the maximum of (24.5) by choosing $n+rx=x\lambda_n$; thus
$$\Lambda(A)=x\lambda_n={-n\lambda_n\over r-\lambda_n}\,,$$
which is the right-hand side of (24.2). By (6.3) and Theorem~12 this
is an upper bound on~$\vartheta$.\ \pfbox

\medskip
Incidentally, we need to be a little careful in (24.2): The
denominator can be zero, but only when $G=\overline{K}_n$.

\meno
{\bf 25. Automorphisms.}\quad
An automorphism of a graph $G$ is a permutation~$p$ of the vertices
such that
$$p(u)\adj p(v)\qquad\hbox{iff }u\adj v\,.\eqno(25.1)$$
Such permutations are closed under multiplication, so they form a
group.

We call $G$ {\it vertex-symmetric\/} if its automorphism group is
vertex-transitive, i.e., if given $u$ and~$v$ there is an
automorphism~$p$ such that $p(u)=v$. We call $G$ {\it
edge-symmetric\/} if its automorphism group is edge-transitive, i.e.,
if given $u\adj v$ and $u'\adj v'$ there is an automorphism~$p$ such
that $p(u)=u'$ and $p(v)=v'$ or $p(u)=v'$ and $p(v)=u'$.

Any vertex-symmetric graph is regular, but edge-symmetric graphs
need not be regular. For example,
$$\unitlength=10pt
\def\putdisk(#1,#2){\put(#1,#2){\disk{.4}}}
\eqalign{\vcenter
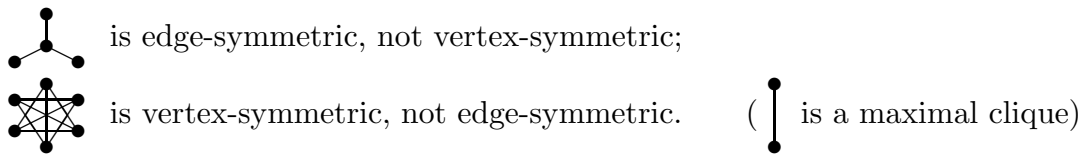
{
\hbox{\beginpicture(2,1.5)(0,0)
\putdisk(0,0)
\putdisk(2,0)
\putdisk(1,.5)
\putdisk(1,1.5)
\put(0,0){\line(2,1){1}}
\put(2,0){\line(-2,1){1}}
\put(1,.5){\line(0,1){1}}
\endpicture}}
\quad&\hbox{is edge-symmetric, not vertex-symmetric;}\cr
\noalign{\smallskip}
\vcenter{
\hbox{\beginpicture(2,2)(0,0)
\putdisk(1,0)
\putdisk(1,2)
\putdisk(0,.5)
\putdisk(0,1.5)
\putdisk(2,.5)
\putdisk(2,1.5)
\put(0,.5){\line(2,1){2}}
\put(2,.5){\line(-2,1){2}}
\put(0,.5){\line(2,3){1}}
\put(2,.5){\line(-2,3){1}}
\put(0,.5){\line(1,0){2}}
\put(0,1.5){\line(1,0){2}}
\put(1,0){\line(-2,3){1}}
\put(1,0){\line(2,3){1}}
\put(1,0){\line(0,1){2}}
\endpicture}}
\quad&\hbox{is vertex-symmetric, not edge-symmetric.}\qquad
 (\vcenter{\hbox{\beginpicture(1,2)(0,0)
\putdisk(.5,0)\putdisk(.5,2)\put(.5,0){\line(0,1){2}}\endpicture}}
\hbox{ is a maximal clique})\cr}$$

\vskip-3pt\noindent
The graph
$\overline{C}_n$ is not edge-symmetric for $n>7$ because it has more edges
than automorphisms. Also, $\overline{C}_7$ has no automorphism that
takes $0\adj 2$ into $0\adj 3$.

\proclaim
Lemma. If $G$ is edge-symmetric and regular, equality holds in
Lemma~24.

\noindent
{\bf Proof.}\quad
Say that $A$ is an optimum feasible matrix for~$G$ if it is a feasible
matrix with
$$\Lambda(A)=\vartheta(G)$$
as in section~6. 
We can prove that optimum feasible matrices form a convex set, as follows.
First, $tA+(1-t)B$ is clearly feasible when $A$ and~$B$ are feasible.
Second,
$$\Lambda\bigl(tA+(1-t)B\bigr)\leq t\Lambda(A)+(1-t)\Lambda(B)\,,\qquad
0\leq t\leq 1\eqno(25.2)$$
holds for all symmetric matrices $A$ and~$B$, by (6.2); this follows because
there is a unit vector~$x$ such that
$\Lambda\bigl(tA+(1-t)B\bigr)=x^T\bigl(tA+(1-t)B\bigr)x=tx^T\!Ax+(1-t)x^TBx
\leq t\Lambda(A)+(1-t)\Lambda(B)$. Third, if $A$ and~$B$ are optimum feasible
matrices, the right side of (25.2) is $\vartheta(G)$ while the left side is
$\geq\vartheta(G)$ by (6.3). Therefore equality holds.

If $A$ is an optimum feasible matrix for~$G$, so is $p(A)$, the matrix
obtained by permuting rows and columns by an automorphism~$p$. (I~mean
$p(A)_{uv}=A_{p(u)p(v)}$.) Therefore the average, $\bar{A}$, over
all~$p$ is also an optimal feasible matrix. Since $p(\bar{A})=\bar{A}$
for all automorphisms~$p$, and since $G$ is edge-symmetric,
$\bar{A}$~has the form $J+xB$ where $B$ is the adjacency matrix
of~$G$. The bound in Lemma~24 is therefore tight.\ \pfbox

(Note: If $p$ is a permutation, let $P_{uv}=1$ if $u=p(v)$,
otherwise~0. Then $(P^T\!AP)_{uv}=\sum
P^T_{uj}A_{jk}P_{kv}=A_{p(u)p(v)}$, so $p(A)=P^T\!AP$.)

The argument in this proof shows that the set of all optimum feasible
matrices~$A$ for~$G$ has a common eigenvector~$x$ such that 
$Ax=\vartheta(G)x$.  The argument also shows that, if $G$ has an edge
automorphism taking $u\adj v$ into $u'\adj v'$, we can assume without loss of
generality that $A_{uv}=A_{u'v'}$ in an optimum feasible matrix. This simplifies
the computation of $\Lambda(A)$, and justifies our restriction to circulant
matrices (22.14) in the case of cyclic graphs.

\proclaim
Theorem. If $G$ is vertex-symmetric,
$\vartheta(G)@\vartheta(\,\overline{G}\,)=n$.

\noindent
{\bf Proof.}\quad
Say that $b$ is an optimum normalized labeling of $\overline{G}$ if it is
a normalized orthogonal labeling of~$\overline{G}$ achieving equality in
(7.1) when all weights are~1:
$$\vartheta=\sum_{u,v}b_u\cdot b_v\,,\qquad
\sum_v\Vert b_v\Vert^2=1\,,\qquad
b_u\cdot b_v=0\hbox{ when }u\adj v\,.\eqno(25.3)$$
Let $B$ be the corresponding spud; i.e., $B_{uv}=b_u\cdot b_v$ and
$\vartheta=\sum_{u,v}B_{uv}$. Then $p(B)$ is also equivalent to an
optimum normalized labeling, whenever $p$ is an automorphism;
 and such matrices~$B$ form a convex set,
so we can assume as in the lemma that $B=p(B)$ for all
automorphisms~$p$. Since $G$ is vertex-symmetric, we must have
$B_{vv}=1/n$ for all vertices~$v$. Thus there is an optimum normalized
labeling~$b$ with $\Vert b_v\Vert^2=1/n$, and the arguments of
Lemma~10 and Theorem~12 establish the existence of
such a~$b$ with
$$c(b_v)=\vartheta(G)/n\eqno(25.4)$$ 
for all $v$. But $b$ is an orthogonal labeling of~$\overline{G}$, hence
$$\vartheta_1(\,\overline{G},\ones)\leq n/\vartheta(G)$$
by the definition (5.2) of $\vartheta_1$. Thus $\vartheta(\,\overline{G}\,)
\vartheta(G)\leq n$; we have already proved the reverse inequality in
Lemma~23. \ \pfbox

\meno
{\bf 26. Consequence for eigenvalues.}\quad
A curious corollary of the results just proved is the following fact
about eigenvalues.

\proclaim
Corollary. If the graphs $G$ and~$\overline{G}$ are vertex-symmetric and
edge-symmetric, and if the adjacency matrix of~$G$ has eigenvalues
$$\lambda_1\geq \lambda_2\geq \cdots\geq \lambda_n\,,\eqno(26.1)$$
then
$$(\lambda_1-\lambda_n)(n-\lambda_1+\lambda_2)
=-\lambda_n(\lambda_2+1)n\,.\eqno(26.2)$$

\noindent
{\bf Proof.}\quad
By Lemma 25 and Theorem 25,
$${n\Lambda(-B)\over\Lambda(B)+\Lambda(-B)}\
{n\Lambda(-\overline{B}\,)\over \Lambda(\,\overline{B}\,)
+\Lambda(-\overline{B}\,)} =n\,,\eqno(26.3)$$
where $B$ and $\overline{B}$ are the adjacency matrices of $G$
and~$\overline{G}$, and where we interpret 0/0 as~1. We have
$$\overline{B}=J-I-B\,.\eqno(26.4)$$
If the eigenvalues of $B$ are given by (26.1), the eigenvalues
of~$\overline{B}$ are therefore
$$n-1-\lambda_1\geq -1-\lambda_n\geq\cdots\geq
-1-\lambda_2\,.\eqno(26.5)$$ 
(We use the fact that $G$ is regular of degree~$\lambda_1$.) Formula
(26.2) follows if we plug the values $\Lambda(B)=\lambda_1$,
$\Lambda(-B)=-\lambda_n$, $\Lambda(\,\overline{B}\,)=n-1-\lambda_1$,
$\Lambda(-\overline{B}\,)=1+\lambda_2$ into (26.3).\ \pfbox

\meno
{\bf 27. Further examples of symmetric graphs.}\quad
Consider the graph $P(m,t,q)$ whose vertices are all ${m\choose t}$
subsets of cardinality~$t$ of some given set~$S$ of cardinality~$m$,
where
$$u\adj v\ \Longleftrightarrow\ \vert u\cap v\vert=q\,.\eqno(27.1)$$
We want $0\leq q<t$ and $m\geq 2t-q$, so that the graph isn't empty.
In fact, we can assume that $m\geq 2t$, because $P(m,r,q)$ is
isomorphic to $P(m,m-t,m-2t+q)$ if we map each subset~$u$ into the set
difference $S\setminus u$:
$$\vert(S\setminus u)\cap (S\setminus v)\vert = \vert S\vert-\vert
u\cup v\vert =\vert S\vert -\vert u\vert - \vert v\vert +\vert u\cap
v\vert \,.\eqno(27.2)$$
The letter $P$ stands for Petersen, because $P(5,2,0)$ is the well
known ``Petersen graph'' on 10~vertices,

$$
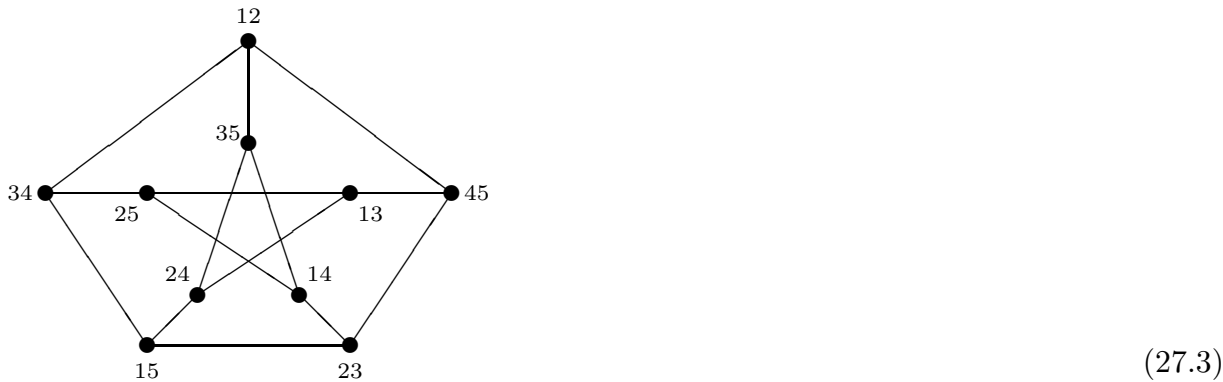
\advance\abovedisplayskip-20pt
\unitlength=16pt
\def\putdisk(#1,#2){\put(#1,#2){\disk{.3}}}
\vcenter{
\beginpicture(8,7)(0,0)
\putdisk(0,3)\putdisk(2,3)\putdisk(6,3)\putdisk(8,3)
\putdisk(4,6)\putdisk(4,4)\putdisk(3,1)\putdisk(5,1)
\putdisk(2,0)\putdisk(6,0)
\sevenrm
\put(-.5,3){\makebox(0,0){34}}
\put(1.6,2.6){\makebox(0,0){25}}
\put(6.4,2.6){\makebox(0,0){13}}
\put(8.5,3){\makebox(0,0){45}}
\put(3.6,4.2){\makebox(0,0){35}}
\put(4,6.5){\makebox(0,0){12}}
\put(2.6,1.4){\makebox(0,0){24}}
\put(5.4,1.4){\makebox(0,0){14}}
\put(2,-.5){\makebox(0,0){15}}
\put(6,-.5){\makebox(0,0){23}}
\put(0,3){\line(1,0){8}}
\put(0,3){\line(4,3){4}}
\put(8,3){\line(-4,3){4}}
\put(4,4){\line(0,1){2}}
\put(3,1){\line(1,3){1}}
\put(5,1){\line(-1,3){1}}
\put(3,1){\line(3,2){3}}
\put(5,1){\line(-3,2){3}}
\put(2,0){\line(1,1){1}}
\put(6,0){\line(-1,1){1}}
\put(2,0){\line(1,0){4}}
\put(2,0){\line(-2,3){2}}
\put(6,0){\line(2,3){2}}
\endpicture}\eqno(27.3)$$

\noindent
These graphs are clearly vertex-symmetric and edge-symmetric, because
every permutation of~$S$ induces an automorphism. For example, to find
an automorphism that maps $u\adj v$ into $u'\adj v'$, let $u=(u\cap
v)\cup\bar{u}$, $v=(u\cap v)\cup \bar{v}$, $u'=(u'\cap
v')\cup\bar{u}'$, $v'=(u'\cap v')\cup \bar{v}'$, and apply any
permutation that takes the $q$~elements of $u\cap v$ into the
$q$~elements of $u'\cap v'$, the $t-q$ elements of~$\bar{u}$ into the
$t-q$ elements of~$u'$, and the $t-q$ elements of $\bar{v}$
into~$\bar{v}'$. Thus we can determine $\vartheta\bigl(P(m,t,q)\bigr)$
from the eigenvalues of the adjacency matrix. 
Lov\'asz [\Lov] discusses the case $q=0$, and his discussion readily
generalizes to other values of~$q$. 
It turns out that $\vartheta\bigl(P(m,t,0)\bigr)={m-1\choose t-1}$. This is
also the value of $\alpha\bigl(P(m,t,0)\bigr)$, because the ${m-1\choose t-1}$
vertices containing any given point form a stable set.

The special case $t=2$, $q=0$ is
especially interesting because those graphs also satisfy the condition
of Corollary~26. We have
$$n={m\choose 2}\,,\qquad \lambda_1={m-2\choose 2}\,,\qquad
\lambda_2=1\,,\qquad \lambda_n=3-m\,,\eqno(27.4)$$
and (26.2) does indeed hold (but not ``trivially'').
It is possible to cover $P(m,2,0)$ with disjoint maximum cliques; hence
$\kappa\bigl(P(m,2,0)\bigr)={m\choose 2}\left/\lfloor{m\over 2}\rfloor\right.
=2\lceil{m\over 2}\rceil -1$. In particular, when $G$~is the Petersen graph we
have $\alpha(G)=\vartheta(G)=4$, $\kappa(G)=5$; also $\alpha(\overline{G})=2$,
$\vartheta(\overline{G})=\kappa(\overline{G})={5\over 2}$.

\meno
{\bf 28. A bound on $\vartheta$.}\quad
The paper [\Lov] contains one more result about $\vartheta$ that is not
in~[\GLSbook], so we will wrap up our discussion of [\Lov] by describing~[\Lov,
Theorem~11].

\proclaim
Theorem. If $G$ has an orthogonal labeling of dimension~$d$ with no
zero vectors, we have $\vartheta(G)\leq d$. 

\noindent
{\bf Proof.}\quad
Given a non-zero orthogonal labeling $a$ of dimension~$d$, we can
assume that $\Vert a_v\Vert^2=1$ for all~$v$. (The hypothesis about
zeros is important, since there is trivially an orthogonal labeling of
any desired dimension if we allow zeros. The labeling needn't be
optimum.) Then we construct an orthogonal labeling $a''$ of
dimension~$d^2$, with $c(a''_v)=1/d$ for all~$v$, as follows:

Let $a'_v$ have $d^2$ components where the $(j,k)$ component is
$a_{jv}a_{kv}$. Then
$$a'_u\cdot a'_v=(a_u\cdot a_v)^2\eqno(28.1)$$
as in (20.7). Let $Q$ be any orthogonal matrix with $d^2$~rows and
columns, such that the $(j,k)$ entry in row $(1,1)$ is $1/\sqrt{d}$ for
$j=k$, 0~otherwise. Then we define
$$a''_v=Qa'_v\,.\eqno(28.2)$$
Once again $a''_u\cdot a''_v=(a_u\cdot a_v)^2$, so $a''$ is an
orthogonal labeling. We also have first component
$$a_{(1,1)v}''=\sum_{j,k}\,{[j=k]\over\sqrt{d}}\;a'_{(j,k)v}
=\sum_k\;{a^2_{kv}\over\sqrt{d}\,}={1\over\sqrt{d}\,}\,;\eqno(28.3)$$
hence $c(a''_v)=1/d$. This proves $\vartheta(G)\leq d$, by definition
of~$\vartheta_1$. \ \pfbox

\medskip
This theorem improves the obvious lower bound $\alpha(G)$ on the dimension of
an optimum labeling.

\meno
{\bf 29. Compatible matrices.}\quad
There's another way to formulate the theory we've been developing,
by looking at things from a somewhat higher 
level, following ideas developed by Lov\'asz and Schrijver [\LS] a few years
after the book~[\GLSbook] was written. Let us say that the matrix $A$ is
{\it$\lambda$-compatible with $G$ and $w$} if $A$ is an $(n+1)\times(n+1)$
spud indexed by the vertices of~$G$ and by a special value~$0$, having the
following properties:
\itemitem{$\bullet$}$A_{00}=\lambda$;
\itemitem{$\bullet$}$A_{vv}=A_{0v}=w_v$ for all vertices~$v$;
\itemitem{$\bullet$}$A_{uv}=0$ whenever $u\noadj v$ in $G$.

\proclaim
Lemma. There exists an orthogonal labeling $a$ for $G$ with costs
$c(a_v)=w_v/\lambda$ if and only if there exists a matrix~$A$ that is
$\lambda$-compatible with $G$ and $w$.

\noindent
{\bf Proof.}\quad
Given such an orthogonal labeling, we can normalize each vector so that
$\Vert a_v\Vert^2=w_v$. Then when $w_v\ne0$ we have
$${w_v\over\lambda}=c(a_v)={a_{1v}^2\over w_v},$$
so we can assume that $a_{1v}=w_v/\sqrt\lambda$ for all~$v$. Add a new
vector~$a_0$, having $a_{10}=\sqrt\lambda$ and $a_{j0}=0$ for all $j>1$.
Then the matrix~$A$ with $A_{uv}=a_u\cdot a_v$ is easily seen to be
$\lambda$-compatible with $G$ and $w$.

Conversely, if such a matrix $A$ exists, there are $n+1$ vectors $a_0$,
\dots,~$a_n$ such that $A_{uv}=a_u\cdot a_v$; in particular,
$\Vert a_{0}\Vert^2=\lambda$. Let $Q$ be an orthogonal matrix with first
row $a_0^T\!/\sqrt\lambda$, and define $a'_v=Qa_v$ for all~$v$. Then
$a'_{10}=\sqrt\lambda$ and $a'_{j0}=0$ for all $j>1$. Also $a'_u\cdot a'_v=
a_u\cdot a_v=A_{uv}$ for all $u$ and~$v$. Hence $\sqrt\lambda a'_{1v}=
a'_0\cdot a'_v=A_{0v}=w_v$ and $\Vert a'_v\Vert^2=a'_v\cdot a'_v=A_{vv}=w_v$,
for all $v\in G$, proving that $c(a'_v)=w_v/\lambda$. Finally $a'$ is
an orthogonal labeling, since $a'_u\cdot a'_v=A_{uv}=0$ whenever
$u\noadj v$.\quad\pfbox

\proclaim
Corollary. $x\in{\tt TH}(G)$ iff there exists a matrix $1$-compatible with
$\overline G$ and~$x$.

\noindent
{\bf Proof.}\quad
Set $\lambda=1$ in the lemma and apply Theorem 14.\quad\pfbox

\medskip
The corollary and definition (4.1) tell us that $\vartheta(G,w)$ is
$\max(w_1x_1+\cdots+w_nx_n)$ over all $x$ that appear in matrices that
are $1$-compatible for $\overline G$ and $x$. Theorem~12 tells us that
$\vartheta(G,w)$ is also the minimum~$\lambda$ such that there exists a
$\lambda$-compatible matrix for $G$ and~$w$. The ``certificate''
property of Theorem~13 has an even stronger formulation in matrix terms:

\proclaim
Theorem. Given a nonnegative weight vector $w=(w_1,\ldots,w_n)^T$, let
$A$ be $\lambda$-compatible with $G$ and~$w$, where $\lambda$ is
as small as possible, and let $B$ be $1$-compatible with $\overline G$
and~$x$, where $w_1x_1+\cdots+w_nx_n$ is as large as possible. Then
$$A@@DB=0\,,\eqno(29.1)$$
where $D$ is the diagonal matrix with $D_{00}=-1$ and $D_{vv}=+1$ for
all $v\ne0$. Conversely, if $A$ is $\lambda$-compatible with $G$ and~$w$
and if $B$ is $1$-compatible with $\overline G$ and~$x$, then (29.1)
implies that $\lambda=w_1x_1+\cdots+w_nx_n=\vartheta(G,w)$.

\noindent
{\bf Proof.}\quad
Assume that $A$ is $\lambda$-compatible with $G$ and $w$, and $B$ is
$1$-compatible with $\overline G$ and~$x$. Let $B'=DBD$, so that $B'$ is
a spud with $B'_{00}=1$, $B'_{0v}=B'_{v0}=-x_v$, and $B'_{uv}=B_{uv}$
when $u$ and~$v$ are nonzero. Then the dot product $A\cdot B'$ is
$$\lambda-w_1x_1-\cdots-w_nx_n-w_1x_1-\cdots-w_nx_n+w_1x_1+\cdots+w_nx_n
 =\lambda-(w_1x_1+\cdots+w_nx_n)\,,$$
because $A_{uv}B_{uv}=0$ when $u$ and $v$ are vertices of $G$. We showed
in the proof of F9 in section~9 that the dot product of spuds is
nonnegative; in fact, that proof implies that the dot product is zero
if and only if the ordinary matrix product is zero. So
$\lambda=w_1x_1+\cdots+w_nx_n=\vartheta(G,w)$ iff $AB'=0$, and this
is equivalent to~(29.1).\quad\pfbox

\medskip
Equation (29.1) gives us further information about the orthogonal labelings
$a$ and~$b$ that appear in Theorems 12 and~13. Normalize those labelings
so that $\Vert a\Vert^2=w_v$ and $\Vert b\Vert^2=x_v$. Then we have
$$\openup 1\jot
\eqalignno{\sum_{t\in G}w_t\,(b_t\cdot b_v)&=\vartheta\,x_v\,,&(29.2)\cr
           \sum_{t\in G}x_t\,(a_t\cdot a_v)&=w_v\,,&(29.3)\cr
           \sum_{t\in G}(a_t\cdot a_u)@(b_t\cdot
           b_v)&=w_ux_v\,,&(29.4)\cr}$$
for all vertices $u$ and $v$ of $G$. (Indeed, (29.2) and (29.4) are
respectively equivalent to $(AB')_{0v}=0$ and $(AB')_{vv}=0$;
(29.3) is equivalent to $(B'A)_{0v}=0$.) Notice that if $\widehat A$ and
$\widehat B$ are the $n\times n$ spuds; obtained by deleting row~0 and
column~0 from optimum matrices $A$ and~$B$, these equations are equivalent to
$$\widehat Bw=\vartheta x\,,\qquad \widehat Ax=w\,,\qquad
\widehat A\widehat B=wx^T\,.\eqno(29.5)$$
Equation (29.1) is equivalent to (29.5) together with the condition
$w_1x_1+\cdots+w_nx_n=\vartheta$.

Since $AB'=0$ iff $B'A=0$ when $A$ and $B'$ are symmetric matrices, the
optimum matrices $A$ and~$B'$ commute. This implies that they have common
eigenvectors: There is an orthogonal matrix~$Q$ such that
$$A=Q\,{\rm diag}\,(\lambda_0,\ldots,\lambda_n)\,Q^T\,,\qquad\qquad
 B'=Q\,{\rm diag}\,(\mu_0,\ldots,\mu_n)\,Q^T\,.\eqno(29.6)$$
Moreover, the product is zero, so
$$\lambda_0\mu_0=\cdots=\lambda_n\mu_n=0\,.\eqno(29.7)$$
The number of zero eigenvalues $\lambda_k$ is $n+1-d$, where $d$ is the
smallest dimension for which there is an orthogonal labeling~$a$ with
$A_{uv}=a_u\cdot a_v$. A similar statement holds for $B'$, since the
eigenvalues of $B$ and $B'$ are the same; $y$ is an eigenvector for $B$ iff
$Dy$ is an eigenvector for $B'$. In the case $G=C_n$, studied in
section~22, we constructed an orthogonal labeling (22.3) with only three
dimensions, so all but~3 of the eigenvalues~$\mu_k$ were zero. When all the
weights~$w_v$ are nonzero and $\vartheta(G)$ is large, Theorem~28 implies
that a large number of $\lambda_k$ must be nonzero, hence a large number of
$\mu_k$ must be zero.

The ``optimum feasible matrices'' $A$ studied in section 6 are related to
the matrices $\widehat A$ of (29.5) by the formula
$$\vartheta \widehat A=ww^T-\vartheta\,{\rm diag}\,(w_1,\ldots,w_n)
-{\rm diag}\,(\sqrt{w_1},\ldots,\sqrt{w_n}\,)\,A\,
 {\rm diag}\,(\sqrt{w_1},\ldots,\sqrt{w_n}\,)\,,\eqno(29.8)$$
because of the construction following (6.6). If the largest eigenvalue
$\Lambda(A)=\vartheta$ of~$A$ occurs with multiplicity~$r$, the rank
of $\vartheta I-A$ will be $n-r$, hence $\widehat A$ will have rank
$n-r$ or $n-r+1$, and the number of zero eigenvalues $\lambda_k$ in
(29.6) will be $r+1$ or $r$.

\meno
{\bf 30. Antiblockers.}\quad
The convex sets {\tt STAB}, {\tt TH}, and {\tt QSTAB} defined in
section~2 have many special properties. For example, they are always
 nonempty, closed, convex, and nonnegative; they also satisfy the condition
$$0\leq y\leq x \quad{\rm and}\quad x\in X\ \Rightarrow\ y\in
X\,.\eqno(30.1)$$ 
A set $X$ of vectors satisfying all five of these properties is called a {\it
convex corner}.

If $X$ is any set of nonnegative vectors we define its {\it
antiblocker\/} by the condition
$${\rm abl}\,X=\{\,y\geq 0\mid x\cdot y\leq 1\hbox{ for all }x\in
X\,\}\,.\eqno(30.2)$$
Clearly abl$\,X$ is a convex corner, and ${\rm abl}\,X\supseteq{\rm abl}\,X'$
when $X\subseteq X'$. 

\proclaim
Lemma. If $X$ is a convex corner we have ${\rm abl}\,{\rm abl}\,X=X$.

\noindent
{\bf Proof.}\quad
(Compare with the proof of F5 in section~8.) The relation
$X\subseteq {\rm abl}\,{\rm abl}\,X$ is obvious by definition (30.2),
so the lemma can fail only if there is some $z\in {\rm abl}\,{\rm
abl}\,X$ with $z\notin X$. Then there is a hyperplane separating~$z$
from~$X$, by~F1; i.e., there is a vector~$y$ and a number~$b$
such that $x\cdot y\leq b$ for all $x\in X$ but $z\cdot y>b$. Let $y'$
be the same as~$y$ but with all negative components changed to zero.
Then $(y',b)$ is also a separating hyperplane. [{\it Proof:\/} If $x\in
X$, let $x'$ be the same as~$x$ but with all components changed to
zero where $y$ has a negative entry; then $x'\in X$, and $x\cdot
y'=x'\cdot y\leq b$. Furthermore $z\cdot y'\geq z\cdot y>b$.] If
$b=0$, we have $\lambda y'\in {\rm abl}\,X$ for all $\lambda >0@$; this
contradicts $z\cdot\lambda y'\leq 1$. 
We cannot have $b<0$, since $0\in X$.
Hence $b>0$, and the vector
$y'/b\in {\rm abl}\,X$. But then $z\cdot (y'/b)$ must be $\leq 1$,
a~contradiction. \ \pfbox

\proclaim
Corollary. If $G$ is any graph we have
$${\tt STAB}(\overline{G})={\rm abl}\;{\tt QSTAB}(G)\,,\eqno(30.3)$$
\vskip-20pt
$${\tt TH}(\overline{G})={\rm abl}\;{\tt TH}(G)\,,\eqno(30.4)$$
\vskip-20pt
$${\tt QSTAB}(\overline{G})={\rm abl}\;{\tt STAB}(G)\,.\eqno(30.5)$$

\noindent
{\bf Proof.}\quad
First we show that
$${\rm abl}\,X=\hbox{abl convex hull $X$}.\eqno(30.6)$$
The left side surely contains the right. And any element $y\in{\rm
abl}\,X$ will satisfy
$$(\alpha_1x^{(1)}+\cdots +\alpha_kx^{(k)})\cdot y\leq 1$$
when the $\alpha$'s are nonnegative scalars summing to~1 and the
$x^{(j)}$ are in~$X$. This proves (30.6), because the convex hull
of~$X$ is the set of all such $\alpha_1x^{(1)}+\cdots
+\alpha_kx^{(k)}$. 

Now (30.6) implies (30.5), because the definitions in section~2 say
that
$$\eqalign{{\tt QSTAB}(\overline{G})&={\rm abl}\,\{\,x\mid x\hbox{ is a
clique labeling of $\overline{G}$}\,\}\cr
\noalign{\smallskip}
&={\rm abl}\,\{\,x\mid x\hbox{ is a stable labeling of
$G$}\,\}\,,\cr
\noalign{\smallskip}
{\tt STAB}(G)&=\hbox{convex hull}\,\{\,x\mid x\hbox{ is a stable
labeling of $G$}\,\}\,.\cr}$$
And (30.5) is equivalent to (30.3) by the lemma, because ${\tt
STAB}(G)$ is a convex corner. (We must prove (30.1), and it suffices to do
this when $y$ equals~$x$ in all but one component; and in fact by
convexity we may assume that $y$ is~0 in that component; and then we
can easily prove~it, because any subset of a stable set is stable.)

Finally, (30.4) is equivalent to Theorem~14, because ${\tt TH}(G)={\rm
abl}\,\{\,x\mid x_v=c(a_v)$  for some orthogonal labeling 
of~$G\,\}$. \ \pfbox

\medskip
The sets {\tt STAB} and {\tt QSTAB} are polytopes, i.e., they are bounded and
can be defined by a finite number of inequalities. But the antiblocker concept
applies also to sets with curved boundaries. For example, let
$$X=\{\,x\geq 0\bigm| \,\Vert x\Vert\leq 1\,\}\eqno(30.7)$$
be the intersection of the unit ball and the nonnegative orthant.
Cauchy's inequality
implies that $x\cdot y\leq 1$ whenever $\Vert x\Vert \leq 1$ and $\Vert
y\Vert\leq 1$, hence $X\subseteq{\rm abl}\,X$. And if $y\in{\rm abl}\,X$ we
have $y\in X$, since $y\neq 0$ implies $\Vert y\Vert =y\cdot(y/\,\Vert
y\Vert)\leq 1$. Therefore $X={\rm abl}\,X$.

In fact, the set $X$ in (30.7) is the only set that equals its own antiblocker.
If $Y={\rm abl}\,Y$ and $y\in Y$ we have $y\cdot y\leq 1$, hence $Y\subseteq
X$; this implies ${\rm abl}\,Y\supseteq X$.

\meno
{\bf 31. Perfect graphs.}\quad
Let $\omega(G)$ be the size of a largest clique in~$G$. The graph~$G$
is called {\it perfect\/} if every induced subgraph~$G'$ of~$G$ can be
colored with $\omega(G')$ colors. (See section~15 for the notion of
induced subgraph. This definition of perfection was introduced by
Claude Berge in 1961.)

Let $G^+$ be $G$ with vertex~$v$ duplicated, as described in
section~16. This means we add
a new vertex~$v'$ with the same neighbors as~$v$ and with $v\adj v'$.

\proclaim
Lemma. If $G$ is perfect, so is $G^+$.

\noindent
{\bf Proof.}\quad
Any induced subgraph of $G^+$ that is not $G^+$ itself is either an induced
subgraph of~$G$ (if it omits $v$ or~$v'$ or both), or has the form $G'^+$ for
some induced subgraph~$G'$ of~$G$ (if it retains $v$ and~$v'$). Therefore it
suffices to prove that $G^+$ can be colored with $\omega(G^+)$ colors.

Color $G$ with $\omega(G)$ colors and suppose $v$ is red. Let $G'$ be the
subgraph induced from~$G$ by leaving out all red vertices except~$v$. Recolor
$G'$ with $\omega(G')$ colors, and assign a new color to the set
$G^+\backslash G'$, which is stable in~$G^+$. This colors $G^+$ with
$\omega(G')+1$ colors, hence $\omega(G^+)\leq\omega(G')+1$.

We complete the proof by showing that $\omega(G^+)=\omega(G')+1$. Let $Q$ be a
clique of size $\omega(G')$ in~$G'$.

{\bf Case 1.} $v\in Q$. Then $Q\cup\{v'\}$ is a clique
of~$G^+$.  

{\bf Case 2.} $v\notin Q$. Then $Q$ contains no red element.

\noindent
In both cases we can conclude that $\omega(G^+)\geq \omega(G')+1$. \ \pfbox

\proclaim
Theorem. If $G$ is perfect, ${\tt STAB}(G)={\tt QSTAB}(G)$.

\noindent
{\bf Proof.}\quad
It suffices to prove that every $x\in{\tt QSTAB}(G)$ with {\it
rational\/} coordinates is a member of ${\tt STAB}(G)$, because ${\tt
STAB}(G)$ is a closed set.

Suppose $x\in{\tt QSTAB}(G)$ and $qx$ has integer coordinates. Let
$G^+$ be the graph obtained from~$G$ by repeatedly duplicating
vertices until each original vertex~$v$ of~$G$ has been replaced by a
clique of size~$qx_v$. Call the vertices of that clique the {\it
clones\/} of~$v$. 

By definition of ${\tt QSTAB}(G)$, if $Q$ is any clique of~$G$ we have
$$\sum_{v\in Q}x_v\leq 1\,.$$
Every clique $Q'$ of $G^+$ is contained in a clique of size
$\sum_{v\in Q}qx_v$ for some clique~$Q$ of~$G$. (Including all clones
of each element yields this possibly larger clique.) Thus
$\omega(G^+)\leq q$, and the lemma tells us that $G^+$ can be colored
with $q$~colors because $G^+$ is perfect.

For each color $k$, where $1\leq k\leq q$, let $x^{(k)}_v=1$ if some
clone of~$v$ is colored~$k$, otherwise $x_v^{(k)}=0$. Then $x^{(k)}$
is a stable labeling. Hence
$${1\over q}\ \sum_{k=1}^qx^{(k)}\in{\tt STAB}(G)\,.$$
But every vertex of $G^+$ is colored, so $\sum_{k=1}^q x_v^{(k)}=qx_v$
for all~$v$, so $q^{-1}\sum_{k=1}^qx^{(k)}=x$.\ \pfbox

\meno
{\bf 32. A characterization of perfection.}\quad
The converse of Theorem 31 is also true; but before we prove it we
need another fact about convex polyhedra.

\proclaim
Lemma. Suppose $P$ is the set $\{\,x\geq 0\mid x\cdot z\leq 1\hbox{
for all }z\in Z\,\}={\rm abl}\,Z$ for some finite set~$Z$ and suppose
$y\in {\rm abl}\,P$, i.e.,
$y$ is a nonnegative vector such that $x\cdot y\leq 1$ for all $x\in
P$. Then the set
$$\postdisplaypenalty=10000
Q=\{\,x\in P\mid x\cdot y=1\,\}\eqno(32.1)$$
is contained in the set $\{\,x\mid x\cdot z=1\,\}$ for some $z\in Z$
(unless $Q$ and~$Z$ are both empty). 

\noindent
{\bf Proof.}\quad
This lemma is ``geometrically obvious''---it says that every vertex,
edge, etc., of a convex polyhedron is contained in some
``facet''---but we ought also to prove~it. The proof is by induction
on~$\vert Z\vert$. If $Z$ is empty, the result holds because $P$ is
the set of all nonnegative~$x$, hence $y$ must be~0 and $Q$ must be
empty.

Suppose $z$ is an element of~$Z$ that does not satisfy the condition;
i.e., there is an element $x\in P$ with $x\cdot y=1$ and $x\cdot z\neq
1$. Then $x\cdot z<1$. Let $Z'=Z\setminus\{z\}$ and $P'={\rm
abl}\,Z'$. It follows that $x'\cdot y\leq 1$ for all $x'\in P'$. For if
$x'\cdot y>1$, a~convex combination $x''=\epsilon x+ (1-\epsilon)x'$
will lie in~$P$ for sufficiently small~$\epsilon$, but $x''\cdot y>1$.

Therefore by induction, $Q'=\{\,x\in P'\mid x\cdot y=1\,\}$ is
contained in $\{\,x\mid x\cdot z'=1\,\}$ for some $z'\in Z'$, unless
$Q'$ is empty, when we can take $z'= z$. And $Q\subseteq Q'$,
since $P\subseteq P'$.  \pfbox

\proclaim
Theorem. $G$ is perfect if and only if ${\tt STAB}(G)={\tt QSTAB}(G)$. 

\noindent
{\bf Proof.}\quad
As in section 15, 
let $G\vert U$ be the graph induced from~$G$ by restriction to
vertices~$U$. If $X$ is a set of vectors indexed by~$V$ and if
$U\subseteq V$, let $X\vert U$ be the set of all vectors indexed
by~$U$ that arise from the vectors of~$X$ when we suppress all
components~$x_v$ with $v\notin U$. Then it is clear that
$${\tt QSTAB}(G\vert U)={\tt QSTAB}(G)\vert U\,,\eqno(32.2)$$
because every $x\in{\tt QSTAB}(G\vert U)$ belongs to ${\tt QSTAB}(G)$ if we
set $x_v=0$ for $v\notin U$, and every $x\in{\tt QSTAB}(G)$ satisfies
$\sum_{v\in Q}x_v\leq 1$ for every clique $Q\subseteq U$. Also
$${\tt STAB}(G\vert U)={\tt STAB}(G)\vert U\,,\eqno(32.3)$$
because every stable labeling of $G\vert U$ is a stable labeling of~$G$ if
we extend it with zeros, and every stable labeling of~$G$ is stable
for $G\vert U$ if we ignore components not in~$U$.

Therefore ${\tt STAB}(G)={\tt QSTAB}(G)$ iff ${\tt STAB}(G')={\tt
QSTAB}(G')$ for all induced graphs. By Theorem~31 we need only prove
that ${\tt STAB}(G)={\tt QSTAB}(G)$ implies $G$ can be colored with
$\omega(G)$ colors.

Suppose ${\tt STAB}(G)={\tt QSTAB}(G)$. Then by Corollary~30,
$${\tt STAB}(\overline{G})={\tt QSTAB}(\overline{G})\,.\eqno(32.4)$$
Let $P={\tt STAB}(\overline{G})$, and let $y=\ones/\omega(G)$. Then
$x\cdot y\leq 1$ whenever $x$ is a clique labeling 
 of~$G$, i.e., whenever $x$ is a stable labeling of~$\overline{G}$;
so $x\cdot y\leq 1$ for all $x\in P$. Let $Z$ be the set of all stable
labelings of~$G$, i.e., clique labelings 
of~$\overline{G}$.  Then $P={\tt
QSTAB}(\overline{G})={\rm abl}\,Z$ and $Z$ is nonempty. So the lemma
applies and it tells us that the set~$Q$ defined in (32.1) is contained in 
 $\{\,x\mid x\cdot z=1\,\}$ for some
stable labeling~$z$ of~$G$. Therefore every maximum clique
labeling~$x$ satisfies $x\cdot z=1$; i.e., every clique of size
$\omega(G)$ intersects the stable set~$S$ corresponding to~$z$. So
$\omega(G')=\omega(G)-1$, where
$$G'=G\vert(V\setminus S)\,.\eqno(32.5)$$
By induction on $\vert V\vert$ we can color the vertices of~$G'$ with
$\omega(G')$ colors, then we can use a new color for the vertices
of~$S$; this colors~$G$ with $\omega(G)$ colors.\ \pfbox

\medskip
Lov\'asz states in [\LL] that he knows no polynomial time algorithm to
test if $G$ is perfect; but he conjectures (``guesses'') that such an
algorithm exists, because the results we are going to discuss next suggest that
much more might be provable.

\meno
{\bf 33. Another definition of $\vartheta$.}\quad
The following result generalizes Lemma~9.3.21 of~[\GLSbook].

\proclaim
Lemma. Let $a$ and $b$ be orthogonal labelings of $G$ and~$\overline{G}$
that satisfy the conditions of Theorem~12, normalized so that
$$\Vert a_v\Vert^2\Vert b_v\Vert^2=w_vc(b_v)\,,\qquad a_{1v}\geq 0\,,\qquad
\hbox{and}\qquad b_{1v}\geq 0\,,\eqno(33.1)$$
for all $v$. Then
$$\sum_va_{jv}b_{kv}=\cases{\sqrt{\vartheta(G,w)}\,,&\rm if $j=k=1$;\cr
0,&\rm otherwise.\cr}\eqno(33.2)$$

\noindent
{\bf Proof.}\quad
Let $a_0=(\sqrt\vartheta,0,\ldots,0)^T$ and $b_0=(-1,0,\ldots,0)^T$. Then
the $(n+1)\times(n+1)$ matrices $A=a^Ta$ and $B=b^Tb$ are spuds, and
$A\cdot B=0$. (In the special case $\Vert a_v\Vert^2=w_v$, 
$\Vert b_v\Vert^2=c(b_v)$, matrix~$B$ is what we called $B'$ in the proof of
Theorem~29.)
Therefore $0={\rm tr}\,A^TB={\rm tr}\,a^Tab^Tb={\rm tr}\,ba^Tab^T=
(ab^T)\cdot(ab^T)$, and we have $ab^T=0$. In other words
$$a_{j0}b_{k0}+\sum_va_{jv}b_{kv}=0$$
for all $j$ and $k$.\quad\pfbox

\medskip
We now can show that $\vartheta$ has yet another definition, in some ways nicer
than the one we considered in section~6. (Someday I~should try to find a
simpler way to derive all these facts.) Call the matrix~$B$ {\it dual
feasible\/} for~$G$ and~$w$ if it is indexed by vertices and
$$\eqalignno{&\hbox{$B$ is real and symmetric;}\cr
&\hbox{$B_{vv}=w_v$ for all $v\in V$;}\cr
&\hbox{$B_{uv}=0$ whenever $u\noadj v$ in $G$;}&(33.3)\cr}$$
and define
$$\eqalignno{
\vartheta_6(G,w)=\max\{\,&\Lambda(B)\,\vert\cr
&B\hbox{ is positive semidefinite and
dual feasible for $G$ and $w$}\,\}\,.&(33.4)\cr}$$%
$\bigl($Compare with the analogous definitions in (6.1) and (6.3).$\bigr)$

\proclaim
Theorem. $\vartheta(G,w)=\vartheta_6(G,w)$.

\noindent{\bf Proof.}\quad
If $B$ is positive semidefinite and dual feasible, and if $\lambda$ is any
eigenvalue of~$B$, we can write $B=QDQ^T$ where $Q$ is orthogonal and $D$ is
diagonal, with $D_{11}=\lambda$. Let $b=\sqrt{D}\,Q^T$; then $b$ is an
orthogonal labeling of~$\overline{G}$ with $\Vert b_v\Vert^2=w_v$ for all~$v$.
Furthermore $c(b_v)=b^2_{1v}/w_v=\lambda\,q^2_{v1}/w_v$, where
$(q_{11},\ldots,q_{n1})$ is the first column of~$Q$. Therefore $\sum_v
c(b_v)w_v=\lambda\sum_v q^2_{v1}=\lambda$, and we have $\lambda \leq
\vartheta_4(G,w)$ by (10.1). This proves that $\vartheta_6\leq \vartheta$. 

Conversely, let $a$ and $b$ be orthogonal labellings of~$G$ and~$\overline{G}$
that satisfy the conditions of Theorem~12. Normalize them so that $\Vert
a_v\Vert^2=c(b_v)$ and $\Vert b_v\Vert^2=w_v$. Then $a^2_{1v}=c(a_v)c(b_v)=
w_vc(b_v)/\vartheta=b^2_{1v}/\vartheta$. The lemma now implies that
$(b_{11},\ldots,b_{1n})^T$ is an eigenvector of~$b^Tb$, with
eigenvalue~$\vartheta$. This proves that $\vartheta\leq \vartheta_6$. \ \pfbox

\proclaim
Corollary. $\vartheta(G)=1+\max\{\,\Lambda(B)/\Lambda(-B)
\,\vert\, B\hbox{ is dual feasible for $G$ and 0}\,\}$.

\noindent{\bf Proof.}\quad
If $B$ is dual feasible for $G$ and~0, its eigenvalues are $\lambda_1\geq
\cdots\geq \lambda_n$ where $\lambda_1=\Lambda(B)$ and
$\lambda_n=-\Lambda(-B)$. Then $B'=I+B/\Lambda(-B)$ has eigenvalues
$1-\lambda_1/\lambda_n,\ldots,\break
1-\lambda_n/\lambda_n=0$. Consequently $B'$ is
positive semidefinite and dual feasible for~$G$ and~$\ones$, and
$1+\Lambda(B)/\Lambda(-B)=\Lambda(B')\leq \vartheta_6(G)$.

Conversely, suppose $B'$ is positive semidefinite and dual feasible for~$G$
and~$\ones$, with $\Lambda(B')=\vartheta=\vartheta(G)$. Let $B=B'-I$. Then $B$
is dual feasible for~$G$ and~0, and $0\leq\Lambda(-B)\leq 1$ since the sum of
the eigenvalues of~$B$ is $\mathop{\rm tr} B=0$. 
It follows that $\vartheta -1=\Lambda(B) \leq\Lambda(B)/\Lambda(-B)$. \ \pfbox

\meno
{\bf 34. Facets of {\tt TH}.}\quad
We know that ${\tt TH}(G)$ is a convex corner set in $n$-dimensional space,
so it is natural to ask whether it might have $(n-1)$-dimensional
facets on its nontrivial boundary---for example, a~straight line segment in two
dimensions, or a region of a plane in three dimensions. This means it
has $n$~linearly independent vectors~$x^{(k)}$ such that
$$\sum_vx_v^{(k)}c(a_v)=1\eqno(34.1)$$
for some orthogonal labeling $a$ of~$G$.

\proclaim
Theorem. If ${\tt TH}(G)$ contains linearly independent solutions
$x^{(1)},\ldots,x^{(n)}$ of (34.1), then there is a 
maximal clique~$Q$ of~$G$ such that
$$c(a_v)=\cases{1\,,&$v\in Q$;\cr
0\,,&$v\notin Q$.\cr}\eqno(34.2)$$

\noindent
{\bf Proof.}\quad
Theorem 14 tell us that every $x^{(k)}\in {\tt TH}(G)$ has
$x_v^{(k)}=c(b_v^{(k)})$
 for some orthogonal labeling of~$\overline{G}$. Set $w_v=
c(a_v)$; then $\vartheta(G,w)=1$, by Theorem 13. We can normalize the
labelings so that $\Vert a_v\Vert^2=a_{1v}=w_v$ and $\Vert b_v^{(k)}\Vert^2
=b_{1v}^{(k)}=x_v^{(k)}$. Hence, by Lemma~33,
$$\sum_vx_v^{(k)}a_v=\pmatrix{1\cr 0\cr \vdots\cr
0\cr}=e_1\,.\eqno(34.3)$$ 

Let
$$Q=\{\,v\mid a_{1v}\neq 0\,\}=\{\,v\mid c(a_v)\neq
0\,\}\eqno(34.4)$$
and suppose $Q$ has $m$ elements. Then (34.3) is equivalent to the
matrix equation
$$Ax^{(k)}=e_1\eqno(34.5)$$
where $A$ is a $d\times m$ matrix and $x^{(k)}$ has $m$ components
$x_v^{(k)}$, one for each $v\in Q$. By hypothesis there are $m$
linearly independent solutions to (34.5), because there are $n$
linearly independent solutions to (34.3). But then there are $m-1$
linearly independent solutions to $Ax=0$, and it follows that $A$ has
rank~1: Every row of~$A$ must be a multiple of the top row (which is
nonzero). And then (34.5) tells us that all rows but the top row are
zero. We have proved that
$$c(a_v)\neq 0\ \Rightarrow \ c(a_v)=1\,.\eqno(34.6)$$
Therefore if $u$ and $v$ are elements of~$Q$ we have $a_u\cdot a_v\neq
0$, hence $u\adj v$; $Q$~is a clique. 

Moreover, $Q$ is maximal. For if $v\notin Q$ is adjacent to all elements
of~$Q$, there is a $k$ such that $x_v^{(k)}>0$. But the characteristic labeling
of $Q\cup\{v\}$ is an orthogonal labeling~$a'$ such that 
$\sum_ux_u^{(k)}c(a'_u)=1+x_v^{(h)}>1$, 
hence $x^{(k)}\notin {\tt TH}(G)$. \ \pfbox

\medskip
Conversely, it is easy to see that the characteristic labeling of any maximal
clique~$Q$ does have $n$~linearly independent vectors satisfying (34.1), so it
does define a facet. For each vertex~$u$ we let $x_u^{(u)}=1$, and
$x_v^{(u)}=0$ for all $v\neq u$ except for one vertex $v\in Q$ with $v\noadj u$
(when $u\notin Q$). Then $x^{(u)}$ is a stable labeling so it is in
${\tt TH}(G)$. The point of the theorem is that a constraint 
$\sum_vx_vc(a_v)\leq 1$ of {\tt TH}$(G)$ that is not satisfied by all
 $x\in{\tt QSTAB}(G)$ cannot correspond to a facet of {\tt TH}$(G)$.

\proclaim
Corollary.
$$\eqalign{{\tt TH}(G)\hbox{ is a 
polytope}\ &\Longleftrightarrow\ {\tt TH}(G)={\tt
QSTAB}(G)\cr
&\Longleftrightarrow\ {\tt TH}(G)={\tt STAB}(G)
\ \Longleftrightarrow \ G \  \hbox{\rm is
perfect}.\cr}$$

\noindent
{\bf Proof.}\quad
If ${\tt TH}(G)$ is a polytope it is defined by facets as in the
theorem, which are nothing more than the constraints of ${\tt
QSTAB}(G)$; hence ${\tt TH}(G)={\tt QSTAB}(G)$. Also the antiblocker
of a convex corner polytope is a polytope, so ${\tt TH}(\overline{G})$ is a
polytope by (30.4); it must be equal to ${\tt QSTAB}(\overline{G})$. Taking
antiblockers, we have ${\tt TH}(G)={\tt STAB}(G)$ by (30.3). The
converses are easy since {\tt STAB} and {\tt QSTAB} are always
polytopes. The connection to perfection is an immediate consequence of
Theorem~32 and Lemma~2. \ \pfbox

\medskip
We cannot strengthen the corollary to say that $\vartheta(G)=\alpha(G)$ holds
if and only if $\vartheta(G)=\kappa(G)$; the Petersen graph (section~27) is a counterexample.

\meno
{\bf 35. Orthogonal labelings in a perfect graph.}\quad
A perfect graph has
$$\vartheta(G,w)=\alpha(G,w)=\max\{x\cdot w\mid \hbox{$x$ is a stable
labeling of $G$}\}\,,\eqno(35.1)$$
and Theorem 12 tells us there exist orthogonal labelings of $G$ and
$\overline G$ such that (12.2) and (12.3) hold. But it isn't obvious
what those labelings might be; the proof was not constructive.

The problem is to find vectors $a_v$ such that $a_u\cdot a_v=0$ when
$u\noadj v$ and such that (12.2) holds; then it is easy to satisfy (12.3)
by simply letting $b$ be a stable labeling where the maximum occurs
in~(35.1).

The following general construction gives an orthogonal labeling (not
necessarily optimum) in any graph: Let $g(Q)$ be a nonnegative number
for every clique~$Q$, chosen so that
$$\sum_{v\in Q} g(Q)=w_v\,,\qquad\hbox{for all $v$}\,.\eqno(35.2)$$
Furthermore, for each clique $Q$, let
$$a_{Qv}=\cases{\sqrt{g(Q)},&if $v\in Q$;\cr
0,&otherwise.\cr}\eqno(35.3)$$
Then
$$a_u\cdot a_v\;=\;\sum_{\{u,v\}\subseteq Q}g(Q)\,,$$
hence $a_u\cdot a_v=0$ when $u\noadj v$. If we also let
$a_{Q0}=\sqrt{g(Q)}$ for all~$Q$, $a_{00}=0$, we find
$$a_0\cdot a_v=a_v\cdot a_v=\sum_{v\in Q}g(Q)=w_v\,.$$
We have constructed a matrix $A$ that is $\lambda$-compatible with $G$
and~$w$, in the sense of section~29, where
$$\lambda=a_0\cdot a_0=\sum_Q g(Q)\,.\eqno(35.4)$$
An orthogonal labeling with costs $c(a'_v)=w_v/\lambda$ can now be
found as in the proof of Lemma~29.

The duality theorem of linear programming tells us that the minimum of
(35.4) subject to the constraints (35.2) is equal to the maximum value
of $w\cdot x$ over all $x$ with $\sum_{v\in Q}x_v\le1$ for all~$Q$.
When $x$ maximizes $w\cdot x$, we can assume that $x\ge0$, because a negative
$x_v$ can be replaced by~0 without decreasing $w\cdot x$ or violating
a constraint. (Every subset of a clique is a clique.) Thus, we are
maximizing $w\cdot x$ over ${\tt QSTAB}(G)$; the construction in the previous
paragraph allows us to reduce $\lambda$ as low as $\kappa(G,w)$.
But $\kappa(G,w)=\vartheta(G,w)$ in a perfect graph, so this
construction solves our problem, once we have computed $g(Q)$.

The special case of a bipartite graph is especially interesting,
because its cliques have only one or two vertices. Suppose all edges
of~$G$ have the form $u\adj v$ where $u\in U$ and $v\in V$, and consider
the network defined as follows: There is a special source vertex~$s$
connected to all $u\in U$ by a directed arc of capacity~$w_u$,
and a special sink vertex~$t$ connected from all $v\in V$ by a directed arc of
capacity~$w_v$. The edges $u\adj v$ of $G$ are also present, directed 
from~$u$ to~$v$ with
infinite capacity. Any flow from $s$ to~$t$ in this network
defines a suitable function~$g$, if we let
$$\eqalign{g(\{u,v\})&=\hbox{the flow in $u\rar v$}\,,\cr
g(\{u\})&=\hbox{$w_u$ minus the flow in $s\rar u$}\,,\cr
g(\{v\})&=\hbox{$w_v$ minus the flow in $v\rar t$}\,,\cr}$$
for all $u\in U$ and $v\in V$. Let $S$ be a subset of $U\cup V$.
If we cut the edges that connect~$s$ or~$t$ with vertices not
in~$S$, we cut off all paths from $s$ to~$t$ if and only if $S$
is a stable set. The minimum cut (i.e., the minimum sum of capacities
of cut edges) is equal to the maximum flow; and it is also equal to
$$\sum_{u\in U}w_u+\sum_{v\in V}w_v-\max\{w\cdot x\mid
\hbox{$x$ is a stable labeling}\}
=\sum_{u\in U}w_u+\sum_{v\in V}w_v-\alpha(G,w)\,.$$
Thus the value of $\lambda=\sum_Qg(Q)$ is
$\sum_{u\in U}w_u-\{\hbox{flow from $s$}\}
+\sum_{v\in V}w_v-\{\hbox{flow to $t$}\}
                 +\{\hbox{flow in $u\rar v$ arcs}\}=\alpha(G,w)=
\vartheta(G,w)$ as desired.

For general perfect graphs $G$, a solution to (35.4) with $\lambda
=\vartheta(G,w)$ can be found in polynomial time as shown in equation (9.4.6)
of [\GLSbook]. 
However, the methods described in [\GLSbook] are not efficient enough for
practical calculation, even on small graphs. 

\meno
{\bf 36. The smallest non-perfect graph.}\quad
The cyclic graph $C_5$ is of particular interest because it is the
smallest graph that isn't perfect, and the smallest case where the
function $\vartheta(G,w)$ is not completely known.

The discussion following Theorem 34 points out that ${\tt TH}(G)$
always has facets in common with ${\tt QSTAB}(G)$, when those facets
belong also to ${\tt STAB}(G)$. It is not hard to see that ${\tt
QSTAB}(C_5)$ has ten facets, defined by $x_j=0$ and
$x_j+x_{j\bmod5}=1$ for $0\le j<5$; and ${\tt STAB}(C_5)$ has
an additional facet defined by $x_0+x_1+x_2+x_3+x_4=2$. The weighted
functions $\alpha$ and $\kappa$ of section~4 are evaluated by
considering the vertices of {\tt STAB} and {\tt QSTAB}:
$$\eqalignno{
\alpha(C_5,\{w_0,\ldots,w_4\}^T)&=\max(w_0+w_2,w_1+w_3,
 w_2+w_4,w_3+w_0,w_4+w_1)\,;&(36.1)\cr
\kappa(C_5,\{w_0,\ldots,w_4\}^T)&=\max\bigl(\alpha(C_5,\{w_0,\ldots,w_4\}^T),
(w_0+\cdots+w_4)/2\bigr)\,.&(36.2)\cr}$$
Where these functions agree, they tell us also the value of $\vartheta$.

For example, let $f(x)=\vartheta(C_5,\{x,1,1,1,1,1\}^T)$. Relations (36.1) and
(36.2) imply that $f(x)=x+1$ when $x\ge2$. Clearly $f(0)=2$, and section~22
tells us that $f(1)=\sqrt5$. Other values of $f(x)$ are not yet known.
Equation (23.2) gives the lower bound $f(x)^2\ge x^2+4$. Incidentally, the
$a$ vectors
$$\openup1\jot
\pmatrix{\sqrt x\cr 1\cr 0\cr 0\cr 0\cr}\qquad
\pmatrix{\sqrt x\cr 1\cr \sqrt{x+1}\cr \sqrt{(x-2)(x+1)}\cr 0\cr}\qquad
\pmatrix{1\cr -\sqrt x\cr 0\cr 0\cr 0\cr}\qquad
\pmatrix{1\cr -\sqrt x\cr 0\cr 0\cr 0\cr}\qquad
\pmatrix{\sqrt x\cr 1\cr -\sqrt{x+1}\cr 0\cr \sqrt{(x-2)(x+1)}\cr}$$
and $b=(1)\,(0)\,(0)\,(0)\,(0)$ establish $f(x)$ for $x\ge2$ in the fashion of
Theorems 12 and~13.

Let $\phi=(1+\sqrt5\,)/2$ be the golden ratio.
The matrices $A$ and $B'$ of Theorem 29, when $G=C_5$ and $w=\ones$, are
$$\openup2\jot
\eqalign{A&=\pmatrix{\sqrt5&1&1&1&1&1\cr
1&1&\phi-1&0&0&\phi-1\cr
1&\phi-1&1&\phi-1&0&0\cr
1&0&\phi-1&1&\phi-1&0\cr
1&0&0&\phi-1&1&\phi-1\cr
1&\phi-1&0&0&\phi-1&1\cr}\,;\cr
B'&={1\over\sqrt5}\pmatrix{\sqrt5&-1&-1&-1&-1&-1\cr
-1&1&0&\phi-1&\phi-1&0\cr
-1&0&1&0&\phi-1&\phi-1\cr
-1&\phi-1&0&1&0&\phi-1\cr
-1&\phi-1&\phi-1&0&1&0\cr
-1&0&\phi-1&\phi-1&0&1\cr}\,.\cr}$$
They have the common eigenvectors
$$\pmatrix{\sqrt5\,\cr 1\cr 1\cr 1\cr 1\cr 1\cr}\qquad
\pmatrix{\sqrt5\,\cr -1\cr -1\cr -1\cr -1\cr -1\cr}\qquad
\pmatrix{0\cr \phi\cr 1\cr -1\cr -\phi\cr 0\cr}\qquad
\pmatrix{0\cr 0\cr \phi\cr 1\cr -1\cr -\phi\cr}\qquad
\pmatrix{0\cr 1\cr -\phi\cr \phi\cr -1\cr 0\cr}\qquad
\pmatrix{0\cr 0\cr 1\cr -\phi\cr \phi\cr -1\cr}\,,$$
with respective eigenvalues $(\lambda_0,\ldots,\lambda_5)=
(2\sqrt5,0,\sqrt5/\phi,\sqrt5/\phi,0,0)$ and
$(\mu_0,\ldots,\mu_5)=(0,2,0,0,1/\phi,1/\phi)$. (Cf.~(29.6) and (29.7).)

\meno
{\bf 37. Perplexing questions.}\quad
The book [\GLSbook] explains how to compute $\vartheta(G,w)$ with given
tolerance~$\epsilon$, in polynomial time using an ellipsoid method, but that 
method is too slow and
numerically unstable to deal with graphs that have more than 10 or so vertices.
Fortunately, however, new ``interior-point methods'' have been developed
for this purpose,
especially by Alizadeh [\Ala,\Alb], who has computed $\vartheta(G)$ when $G$ has
hundreds of vertices and thousands of edges. He has also shown how to find
large stable sets, as a
byproduct of evaluating $\vartheta(G,w)$ when $w$ has integer coordinates.
Calculations on somewhat smaller cyclically symmetric
graphs  have also been reported
by Overton [\Ov]. Further
computational experience with such programs should prove to be very
interesting.

Solutions to
the following four concrete problems may also help shed light on the
subject: 

{\bf P1.}\quad
Describe ${\tt TH}(C_5)$ geometrically. This upright set is isomorphic to its
own anti\-blocker. (Namely, if $(x_0,x_1,x_2,x_3,x_4)\in{\tt TH}(C_5)$, then so
are its cyclic permutations $(x_1,x_2,x_3,x_4,x_0)$, etc., as well as the
cyclic permutations of $(x_0,x_4,x_3,x_2,x_1)$; ${\tt TH}(\overline{C}_5)$
contains the cyclic permutations of $(x_0,x_2,x_4, x_1,x_3)$ and
$(x_0,x_3,x_1,x_4,x_2)$.) Can the values
$f(x)=\vartheta(C_5,\{x,1,1,1,1\}^T)$,
discussed in section 36, be expressed in closed form when $0<x<2$,
using familiar functions?

{\bf P2.}\quad
What is the probable value of $\vartheta(G,w)$ when $G$ is a random graph on
$n$~vertices, where each of the ${n\choose 2}$ possible edges is independently
present with some fixed probability~$p$? (Juh\'asz [\Ju] has solved this
problem in the case $w=\ones$, showing that $\vartheta(G)/\sqrt{(1-p)n/p}$ lies
between~${1\over 2}$ and~2 with probability approaching~1 as $n\rar\infty$.)

{\bf P3.}\quad 
What is the minimum $d$ for which $G$ almost surely 
has an orthogonal labeling of
dimension~$d$ with no zero vectors, when $G$ is a random graph as in
Problem~P2?  (Theorem~28 and the theorem of Juh\'asz [\Ju] show that $d$ must
be at least of order~$\sqrt{n}$. But Lov\'asz tells me that he suspects the
correct answer is near~$n$. Theorem~29 and its consequences might be helpful
here.)

{\bf P4.}\quad
Is there a constant $c$ such that $\vartheta(G)\leq
c\sqrt{n}\,\alpha(G)$ for all
$n$-vertex graphs~$G$? (This conjecture was suggested by Lov\'asz in a recent
letter. He knows no infinite family of graphs where $\vartheta(G)/\alpha(G)$
grows faster than $O(\sqrt{n}/\log n)$. The latter behavior occurs for random
graphs, which have $\alpha(G)=\log_{1/p}n$ with high probability [\Bol,
Chapter~XI].)

\meno
Another, more general, question is to ask whether it is feasible to study
two- or three-dimensional projections of ${\tt TH}(G)$, and whether they have
combinatorial significance. The
function  $\vartheta(G,w)$ gives just a one-dimensional
glimpse. 

Lov\'asz and Schrijver have recently generalized the topics treated here to a
wide variety of more powerful techniques for studying 0--1 vectors
associated with graphs [\LS]. In particular, one of their methods can be
described as follows: Let us say that a {\it strong orthogonal labeling\/}
is a vector labeling such that $\Vert a_v\Vert^2=c(a_v)$ and $a_u\cdot a_v
\ge0$, also satisfying the relation
$$c(a_u)+c(a_v)+c(a_w)-1\le a_u\cdot a_v+a_v\cdot a_w\le c(a_v)\eqno(37.1)$$
whenever $u\noadj w$. In particular, when $w=v$ this relation implies that
$a_u\cdot a_v=0$, so the labeling is orthogonal in the former sense.

Notice that every stable labeling is a strong orthogonal labeling of $\overline
G$. Let $S$ be a stable set and let $u$ and $w$ be vertices such that $u\adj
w$. If $u$ and $w$ are not in~$S$, condition
(37.1) just says that $0\le c(a_v)\le 1$,
which surely holds. If $u$ is in~$S$, then $w\notin S$ and (37.1) reduces to
$c(a_v)\le c(a_v)\le c(a_v)$; this holds even more surely.

Let
$${\tt TH}_-(G)=\{\,x\mid \hbox{$x_v=c(b_v)$ for some strong
 orthogonal labeling of $\overline G$}\,\}.\eqno(37.2)$$
(This set is called $N_+({\tt FR}(G))$ in [\LS].) We also define
$$\vartheta_-(G,w)=\max\{\,w\cdot x\mid x\in {\tt TH}_-(G)\,\}\,.\eqno(37.3)$$
The argument in the two previous paragraphs implies that
$$ {\tt STAB}(G)\subseteq {\tt TH}_-(G)\subseteq {\tt TH}(G)\,,$$
hence
$$\alpha(G,w)\le \vartheta_-(G,w)\le\vartheta(G,w)\,.\eqno(37.4)$$
The authors of 
[\LS] prove that $\vartheta_-(G,w)$ can be computed in polynomial
time, about as easily as $\vartheta(G,w)$; moreover, it can be a significantly
better approximation to $\alpha(G,w)$. They show, for example, that
${\tt TH}_-(G)={\tt STAB}(G)$ when $G$ is any cyclic graph~$C_n$.
In fact, they prove that if
$x\in {\tt TH}_-(G)$ and if $v_0\adj v_1$, $v_1\adj v_2$, \dots,
$v_{2n}\adj v_0$ is any circuit or multicircuit of $G$, then
$x_{v_0}+x_{v_1}+\cdots+x_{v_{2n}}\le n$.
This suggests additional research problems:

{\bf P5.}\quad
 What is the smallest graph such that ${\tt STAB}(G) \ne{\tt TH}_-(G)$?

{\bf P6.}\quad
What is the probable value of $\vartheta_-(G)$ when $G$ is a random graph as in
Problem~P2?

\medskip
A recent theorem by Arora, Lund, Motwani, Sudan, and Szegedy [\ALMSS] proves
that there is an $\epsilon >0$ such that no polynomial algorithm can compute a
number between $\alpha(G)$ and $n^{\epsilon}\alpha(G)$ for all $n$-vertex
graphs~$G$, unless $P=NP$. 
Therefore it would be surprising if the answer
to~P6 turns out to be that $\vartheta_-(G)$ is, say, $O(\log n)^2$ with
probability $\rar 1$ for random~$G$. Still, this would not be inconsistent with
[\ALMSS], because the graphs for which $\alpha(G)$ is hard to approximate might
be decidedly nonrandom.

Lov\'asz has called my attention to papers by Kashin and Konyagin 
[\KK, \Ko], which prove (in a very disguised form, related to (6.2) and
Theorem~33) that if $G$ has no stable set with 3~elements we have
$$\vartheta(G)\leq 2^{2/3}n^{1/3}\,;\eqno(37.5)$$
moreover, such graphs exist with
$$\vartheta(G)=\Omega(n^{1/3}\!/\sqrt{\log n}\,)\,.\eqno(37.6)$$

Further study of methods like those in [\LS] promises to be exciting indeed.
Lov\'asz has sketched yet another approach in~[\LLL].

\bigskip
\centerline{\bf References}

\bib
[\Ala]
Farid Alizadeh, ``A sublinear-time randomized parallel algorithm for the
maximum clique problem in perfect graphs,'' {\sl ACM-SIAM Symposium on Discrete
Algorithms\/ \bf 2} (1991), 188--194.

\bib
[\Alb]
Farid Alizadeh, ``Interior point methods in semidefinite programming with
applications to combinatorial optimization,'' preprint, 1993,
available from {\tt alizadeh@icsi.\allowbreak
berkeley.edu}.

\bib
[\ALMSS]
Sanjeev Arora, Carsten Lund, Rajeev Motwani, Madhu Sudan, and Mario Szegedy,
``Proof verification and intractability of approximation problems,'' {\sl
Proceedings of the 33rd IEEE Symposium on Foundations of Computer Science\/}
(1992), 14--23.

\bib
[\Bol]
B\'ela Bollob\'as, {\sl Random Graphs\/} (Academic Press, 1985).

\bib
[\GLSell]
Martin Gr\"otschel, L. Lov\'asz, and A. Schrijver,
``The ellipsoid method and its consequences in combinatorial
optimization,'' 
{\sl Combinatorica\/ \bf 1} (1981), 169--197.

\bib
[\GLSrelax]
M. Gr\"otschel, L. Lov\'asz. and A. Schrijver,
``Relaxations of vertex packing,'' {\sl Journal of Combinatorial Theory\/ 
\bf B40} (1986), 330--343.

\bib
[\GLSbook]
Martin Gr\"otschel, L\'aszl\'o Lov\'asz, and Alexander Schrijver, 
{\sl Geometric Algorithms
and Combinatorial Optimization\/} (Springer-Verlag, 1988), \S9.3. 

\bib
[\House]
A. S. Householder, ``Unitary triangularization of a nonsymmetric matrix,''
{\sl Journal of the Association for Computing Machinery\/ \bf5} (1958),
339--342.

\bib
[\Ju]
Ferenc Juh\'asz, ``The asymptotic behaviour of Lov\'asz' $\vartheta$~function
for random graphs,'' {\sl Combinatorica\/ \bf 2} (1982), 153--155.

\bib
[\KK]
B. S. Kashin and S. V. Kon\t\i agin,
``O sistemakh vektorov v Gil'bertovom prostranstve,''
{\sl Trudy Matematicheskogo Instituta imeni V.~A. Steklova\/ \bf 157} (1981),
64--67. English translation, ``On systems of vectors in a Hilbert space,''
{\sl Proceedings of the Steklov Institute of Mathematics\/} (American
Mathematical Society, 1983, issue~3), 67--70.

\bib
[\Ko]
S. V. Kon\t\i agin, ``O sistemakh vektorov v Evclidovom prostranstve i 
odno\u\i\ \'ekstremal'\-no\u\i\ zadache dl\t\i a mnogochlenov,''
{\sl Matematicheskie Zametki\/ \bf 29} (1981), 63--74. English translation,
``Systems of vectors in Euclidean space and an extremal problem for
polynomials,'' {\sl Mathematical Notes of the Academy of Sciences of the USSR\/
\bf 29} (1981), 33--39.

\bib
[\Lov]
L. Lov\'asz, ``On the Shannon capacity of a graph,'' {\sl IEEE
Transactions on Information Theory\/} {\bf IT--25} (1979), 1--7.

\bib
[\LL] 
L. Lov\'asz, An Algorithmic Theory of Numbers, Graphs, and
Convexity. {\sl CBMS Regional Conference Series in Applied
Mathematics\/} (SIAM, 1986), \S3.2.

\bib
[\LLL]
L. Lov\'asz, ``Stable sets and polynomials,'' to appear in {\sl Discrete
Mathematics}. 

\bib
[\LS]
L. Lov\'asz and A. Schrijver, ``Cones of matrices and set functions and 0--1
optimization,'' {\sl SIAM Journal on Optimization\/ \bf 1} (1991), 166--190.

\bib
[\Ov] 
Michael Overton, ``Large-scale optimization of eigenvalues,'' 
{\sl SIAM Journal on Optimization\/ \bf 2} (1992), 88--120.

\bib
[\Padberg]
M. W. Padberg, ``On the facial structure of set packing polyhedra,''
{\sl Mathematical Programming\/ \bf 5} (1973), 199--215.

\bib
[\Shannon]
Claude E.  Shannon, ``The zero error capacity of a channel,''
 {\sl IRE Transactions on Information Theory\/
\bf 2}, 3 (Sep 1956), 8--19.

\bye